\documentclass[11pt,reqno]{amsart}
\usepackage{amsmath,amssymb,mathrsfs}
\usepackage{xcolor}
\colorlet{mdtRed}{red!50!black}
\definecolor{dblue}{rgb}{0,0,.6}
\usepackage[colorlinks]{hyperref}
\hypersetup{colorlinks,linkcolor={red},citecolor={blue},urlcolor={red}}
\usepackage[all]{xy}
\usepackage{tikz,tikz-cd,tkz-graph,enumerate}

\DeclareMathOperator{\Pic}{\textnormal{Pic}}

\DeclareMathOperator{\Id}{{\rm Id}}
\DeclareMathOperator{\Hom}{{\rm Hom}}

\DeclareMathOperator{\Br}{\textnormal{Br}}

\DeclareMathOperator{\psp}{\mathrm{PSp}}

\newcommand{\mf}[1]{\mathfrak{#1}}
\newcommand{\mc}[1]{\mathcal{#1}}
\newcommand{\bb}[1]{\mathbb{#1}}

\renewcommand{\ker}{\mathrm{Ker}}
\newcommand{\coker}{\mathrm{Coker}}

\newtheorem{theorem}{Theorem}[section] 
\newtheorem{lemma}[theorem]{Lemma} 
\newtheorem{proposition}[theorem]{Proposition}
\newtheorem{corollary}[theorem]{Corollary}

\theoremstyle{definition}
\newtheorem{definition}[theorem]{Definition}
\newtheorem{remark}[theorem]{Remark}

\numberwithin{equation}{section}

\begin{document}
	
	\baselineskip=15.5pt 
	
	\title[Brauer group of moduli stack of parabolic $\psp(r,\bb{C})$--bundles]{Brauer group of moduli
		stack of parabolic $\psp(r,\bb{C})$--bundles over a curve}
	
	\author[I. Biswas]{Indranil Biswas}
	
	\address{Department of Mathematics, Shiv Nadar University, NH91, Tehsil
		Dadri, Greater Noida, Uttar Pradesh 201314, India}
	
	\email{indranil.biswas@snu.edu.in, indranil29@gmail.com}
	
	\author[S. Chakraborty]{Sujoy Chakraborty}
	
	\address{Department of Mathematics, 
		Indian Institute of Science Education and Research Tirupati, Andhra Pradesh 517507, India}
	\email{sujoy.cmi@gmail.com}
	
	\author[A. Dey]{Arijit Dey}
	
	\address{Department of Mathematics, Indian Institute of Technology Chennai, Chennai, India}
	
	\email{arijitdey@gmail.com}
	
	
	
	\begin{abstract}
		Take an irreducible smooth complex projective curve $X$ of genus $g$, with $g\,\geq\, 3$. Let $r$ be an 
		even positive integer. We prove that the Brauer group of the moduli stack of stable parabolic 
		$\psp(r,\bb{C})$--bundles on $X$, of full-flag parabolic data along a set of marked points on $X$, 
		coincides with the Brauer group of the smooth locus of the corresponding coarse moduli space of stable 
		parabolic $\psp(r,\bb{C})$--bundles. Under certain conditions on the parabolic types, we also compute the 
		Brauer group of the smooth locus of this coarse moduli space. Similar computations are also done
		for the case of partial flags.
	\end{abstract}
	
	\maketitle
	
	\section{Introduction}
	
	The cohomological Brauer group of a quasi-projective variety $Y$ over $\bb{C}$, denoted by $\Br(Y)$, is 
	defined to be the torsion part $H^2_{\text{\'et}}(Y,\,\bb{G}_m)_{tor}\, \subset\,
	H^2_{\text{\'et}}(Y,\,\bb{G}_m)$. When $Y$ is smooth, it is known that 
	$H^2_{\text{\'et}}(Y,\,\bb{G}_m)$ is already torsion. Brauer groups are interesting objects to study for a 
	number of reasons. It is a stable birational invariant for smooth projective varieties defined over a field, 
	making it very useful in studying rationality questions. In fact, Brauer group has been used in constructing 
	examples of non-rational varieties by many, including Colliot-Th\'el$\grave{e}$ne, Saltman, Peyre and 
	others. It also plays a central role in the Brauer--Manin obstruction theory, which deals with the study of 
	rational points on varieties defined over number fields. For an algebraic stack, its Brauer group shall mean
	the cohomological Brauer group.
	
	The study of Brauer groups in the context of moduli of parabolic vector bundles over curves has been carried 
	out in recent times \cite{B24, BB23,BCD23,BCD24,BD11}. The computation of the Brauer group of 
	the moduli space of stable parabolic principal bundles for structure groups $\text{SL}(r,\bb{C})$, 
	$\textnormal{PGL}(r,\bb{C})$ and $\textnormal{Sp}(r,\bb{C})$ have been carried out earlier in \cite{BD11, 
		BCD23, BCD24}. Our aim here is to address the case of $\text{PSp}(r,\bb{C})$.
	
	The set-up is as follows. Let $X$ be an irreducible smooth projective curve over $\bb{C}$ of genus $g$, 
	with $g\,\geq\, 3$. A parabolic vector bundle on $X$, denoted by $E_*$, is an algebraic vector bundle $E$ on $X$ 
	together with the data of weighted filtrations on the fibers of $E$ over finitely many fixed marked points on $X$. 
	A symplectic form on $E_*$ with values in a line bundle $L$ is a non-degenerate skew-symmetric parabolic 
	bilinear form $E_*\otimes E_* \,\longrightarrow\, L$, where $L$ is considered as a parabolic bundle with the trivial 
	parabolic structure (no nonzero parabolic weight is assigned to $L$ at the marked points; see 
	\S~\ref{section:preliminaries} for the details). Consider the moduli space of stable parabolic symplectic vector 
	bundles on $X$ of a fixed rank $r$, where $r$ is an even positive integer. The group of $2$--torsion line 
	bundles on $X$ acts on this moduli space through the operation of tensor product 
	(see \S~\ref{section:fixed-point-locus} for more 
	details). The coarse moduli space of stable parabolic $\psp(r,\bb{C})$--bundles on $X$ is 
	the quotient of the stable parabolic symplectic moduli space under this action of the group of $2$--torsion line
	bundles on $X$.
	
	More precisely, fix an even integer $r\,\geq\, 2$, a finite subset $S\, \subset\, X$ of parabolic points, 
	and a line bundle $L$ on $X$. Fix a system of multiplicities $\boldsymbol{m}$ and a system of weights 
	$\boldsymbol{\alpha}$ at the points of $S$ (see \S~ \ref{section:preliminaries} for details). For $d\,\in\,
	\{0,\,1\}$, let $\mf{N}^{\boldsymbol{m,\alpha},d}_{L}$ denote the moduli stack of stable parabolic 
	$\psp(r,\bb{C})$--bundles on $X$ having system of weights $\boldsymbol{\alpha}$, multiplicities 
	$\boldsymbol{m}$ and topological type $d$; the symplectic form takes values in the line bundle $L$. Let
	$\mc{N}^{\boldsymbol{m,\alpha},d}_L$ denote the corresponding coarse moduli space of stable
	parabolic $\psp(r,\bb{C})$--bundles (see Definition \ref{def:psp-moduli} for more details). Denote 
	$L(S) \,:=\, L\otimes \mc{O}_X(S)$. Our main results are the following:
	
	\begin{theorem}[{See Theorem \ref{thm:brauer-group-of-moduli-stack}}]\label{thmi1}
		Let $\left(\mc{N}^{\boldsymbol{m,\alpha},d}_L\right)^{sm}$ denote the smooth locus of
		$\mc{N}^{\boldsymbol{m,\alpha},d}_L$. When $\boldsymbol{m}$ is a full-flag system of multiplicities,
		\begin{align*}
			\Br\left(\mf{N}^{\boldsymbol{m,\alpha},d}_L\right)\,\ \simeq \,\
			\Br\left(\left(\mc{N}^{\boldsymbol{m,\alpha},d}_L\right)^{sm}\right).
		\end{align*}
	\end{theorem} 
	
	\begin{theorem}[{See Theorem \ref{thm:brauer-group-concentrated-weights}
			and Corollary \ref{cor:arbitrary-weight}}]\label{thmi2}
		Assume that $\boldsymbol{\alpha}$ is a generic system of weights that does not contain $0$. The Brauer
		group of the smooth locus $\left(\mc{N}^{\boldsymbol{m,\alpha},d}_{L(S)}\right)^{sm}\ \subset\
		\mc{N}^{\boldsymbol{m,\alpha},d}_{L(S)}$ has the following description:
		\begin{enumerate}[(1)]
			\item If $d\,=\,0$ (equivalently, $\deg(L)$ is even),\, $\frac{r}{2}\,\geq\, 3$ is odd and $m^i_{_{p}}\,=\, 1$ for
			some $p\,\in\, S$ and some $i$,
			\begin{align}
				\Br\left(\left(\mc{N}^{\boldsymbol{m,\alpha},d}_{L(S)}\right)^{sm}\right)
				\,\ \stackrel{\simeq}{\longrightarrow}\ \, \dfrac{H^2(\Gamma,\,\bb{C}^*)}{\frac{\bb{Z}}{2\bb{Z}}}.
			\end{align}
			
			\item If $d\,=\,0$ (equivalently, $\deg(L)$ is even),\, $\frac{r}{2}\,\geq\, 3$ is even and $m^i_{_{p}}\,=\,1$ for
			some $p\,\in\, S$ and some $i$,
			\begin{align}
				\Br\left(\left(\mc{N}^{\boldsymbol{m,\alpha},d}_{L(S)}\right)^{sm}\right)\ \,\stackrel{\simeq}{\longrightarrow}
				\ \, H^2(\Gamma,\,\bb{C}^*).
			\end{align}
			
			\item If $d\,=\,1$ (equivalently $\deg(L)$ is odd),\, $\frac{r}{2}\,\geq\, 3$ is even and $m^i_{_{p}}\,=\, 1$ for
			some $p\,\in\, S$ and some $i$,
			\begin{align}
				\Br\left(\left(\mc{N}^{\boldsymbol{m,\alpha},d}_{L(S)}\right)^{sm}\right)\ \,\stackrel{\simeq}{\longrightarrow}
				\ \, H^2(\Gamma,\,\bb{C}^*).
			\end{align}
			
			\item If $d\,=\,1$ (equivalently, $\deg(L)$ is odd) and $\frac{r}{2}\,\geq\, 3$ is odd,
			\begin{align}
				\Br\left(\left(\mc{N}^{\boldsymbol{m,\alpha},d}_{L(S)}\right)^{sm}\right)\ \,\stackrel{\simeq}{\longrightarrow}
				\ \, H^2(\Gamma,\,\bb{C}^*).
			\end{align}
		\end{enumerate}
	\end{theorem}
	
	Theorem \ref{thmi1} is proved by obtaining a codimension estimation of the fixed point locus, for the action of
	the group of $2$--torsion line 
	bundles on $X$, on the moduli space of stable parabolic symplectic vector bundles, and then using purity 
	results for Brauer groups as in \cite{Ce19}. Theorem \ref{thmi2} is first proved for a \textit{concentrated} 
	system of weights, and later it is extended to arbitrary generic systems of 
	weights using the existence of certain birational maps between moduli spaces of different systems of weights 
	arising through wall-crossing arguments for variations of parabolic weights \cite{Th96}, \cite{DolHu98}.
	
	\section{the set-up}\label{section:preliminaries}
	
	Fix a smooth irreducible complex projective curve $X$ of genus $g$, with $g\,\geq\, 3$. Fix a
	finite subset $S\, \subset\, X$ of distinct closed
	points; these points are referred to as ``parabolic points''.
	
	\begin{definition}\label{def:parabolic-bundles}
		A \textit{parabolic vector bundle} of rank $r$ on $X$ is an algebraic vector bundle $E$ of
		rank $r$ over $X$ together with the data of a weighted flag on the fiber of $E$ over each $p\,\in\, S$:
		\begin{align}\label{eqn:parabolic-data}
			E_{_p} \,=\, E^1_{_{p}}\,\supsetneq\, E^2_{_{p}}\,\supsetneq\, \cdots \,\supsetneq\, E^{\ell(p)}_{_{p}}
			\,\supsetneq\, E^{\ell(p)+1}_{_{p}} \,=\,0\\
			0\,\leq\, \alpha^1_{_{p}}\,<\,\alpha^2_{_{p}} \,<\,\cdots\,<\,\alpha^{\ell(p)}_{_{p}}\,<\,1,\nonumber
		\end{align}
		where $\alpha^i_{_{p}}$ are real numbers. 
		\begin{itemize}
			\item Such a flag is said to be of length $\ell(p)$, and the numbers $m^i_{_{p}} \,:=\, \dim E^i_{_{p}} -
			\dim E^{i+1}_{_{p}}$ are called the \textit{multiplicities} of the flag at $p$. More precisely,
			$m^i_{_{p}}$ is the multiplicity of the weight $\alpha^i_{_{p}}$.
			
			\item The flag at $p$ is said to be \textit{full} if $m^i_{_{p}} \,=\,1$ for every $i$, in 
			which case clearly we have $\ell(p) \,=\, r$.
			
			\item The collection of real numbers $\boldsymbol{\alpha}\,:=\,\{(\alpha^1_{_{p}}
			\,<\,\alpha^2_{_{p}}\,<\,\cdots\,<\,\alpha^{\ell(p)}_{_{p}})\}_{p\in S}$ is called a system
			of weights. 
			
			\item A \textit{parabolic data} consists of a collection
			$\{(E^{\bullet}_{_{p}},\,\alpha^{\bullet}_{_{p}})\}_{p\in S}$ of weighted flags as above. 
			
			\item Sometimes a system of multiplicities (respectively, a system of weights) will be denoted by the bold symbol 
			$\boldsymbol{m}$ (respectively, $\boldsymbol{\alpha}$), when there is no scope of any confusion. Also, 
			a parabolic vector bundle will often be denoted simply by $E_*$, suppressing the parabolic data.
		\end{itemize}
	\end{definition}
	
	\begin{remark}\label{rem:special-structure}
		Let $E_*$ be a parabolic vector bundle of rank $r$ having the trivial weighted flag at each
		$p\,\in\, S$, i.e., $\ell(p) \,=\, 1$ (so that $E^2_{_{p}} \,=\,0$ in
		\eqref{eqn:parabolic-data}) and $\alpha^1_{_{p}} \,=\,0$ is the single weight at each
		$p\,\in\, S$. In this case, it is said that $E_*$ has the \textit{trivial} parabolic structure.
		We shall not distinguish between a vector bundle $E$ and the parabolic bundle $E$
		equipped with the trivial parabolic structure. 
	\end{remark}
	
	\begin{definition}\label{def:parabolic-morphism}
		Let $E_*$ and $F_*$ be two parabolic vector bundles with systems of multiplicities and
		weights being $(\boldsymbol{m,\,\alpha})$ and $(\boldsymbol{m',\,\alpha'})$
		respectively. A \textit{parabolic morphism} $f_* \ :\ E_* \ \longrightarrow\ F_*$
		is an ${\mathcal O}_X$--linear homomorphism $f \,:\, E\,\longrightarrow\, F$
		between the underlying vector bundles such that for each parabolic point $p$,
		$$\{\alpha^i_{_{p}}\,> \,{\alpha'}^j_{_{p}}\} \,\,\implies\,\,
		\{f_{_{p}}(E^i_{_{p}})\,\subset\, {F}^{j+1}_{_{p}}\}.$$
	\end{definition}
	
	Recall the notion of parabolic symplectic vector bundles on curves following \cite{BMW11} (see also 
	\cite{CM24} for higher dimensions). Fix an even positive integer $r$ and a line bundle $L$ on $X$. Equip $L$ with 
	the trivial parabolic structure as described in Remark \ref{rem:special-structure}. Take a parabolic vector bundle $E_*$
	on $X$ of rank $r$. Let
	$$
	\varphi_*\ :\ E_*\otimes E_*\ \longrightarrow\ L
	$$
	be a skew-symmetric homomorphism of parabolic vector bundles.
	Note that $\mc{O}_X\,=\, \mc{O}_X\cdot {\rm Id}_E\, \subset\,
	E_*\otimes E_*^{\vee}$, where $E^\vee_*$ is the parabolic dual of $E_*$, is a parabolic subbundle with
	the trivial parabolic structure. Let 
	\begin{equation}\label{e2}
		\widehat{\varphi}_* \,\,:\,\, E_*\,\,\longrightarrow\,\, E_*^{\vee}\otimes L\,=\, L\otimes E^\vee_*
	\end{equation}
	be the parabolic morphism defined by the following composition of maps:
	\begin{align*}
		E_* \,\simeq\, E_*\otimes \mc{O}_X\,\hookrightarrow\, E_*\otimes (E_*\otimes E_*^{\vee}) \,=\,
		(E_*\otimes E_*)\otimes E_*^{\vee} \,\xrightarrow{\,\,\, \varphi_*\otimes \text{Id}\,\,}\, L\otimes E_*^{\vee}.
	\end{align*}
	
	\begin{definition}\label{def:parabolic-symplectic-bundle}
		A \textit{parabolic symplectic vector bundle} on $X$ taking values in $L$ is a pair $(E_*,\,\varphi_*)$ as
		above such that $\widehat{\varphi}_*$ in \eqref{e2} is an isomorphism of parabolic vector bundles.
	\end{definition}
	
	\begin{definition}\label{def:parabolic-symplectic-isomorphism}
		Two parabolic symplectic vector bundles $(E_*,\, \varphi_*)$ and $(E'_* ,\, \varphi'_*)$ taking values
		in the same line bundle $L$ are said to be \textit{isomorphic} if there exists an isomorphism
		of parabolic vector bundles $\theta_* \,:\, E_*\, \stackrel{\simeq}{\longrightarrow}\, E'_*$ (see
		Definition \ref{def:parabolic-morphism}) satisfying the condition that the following diagram is commutative:
		\begin{align}\label{cd1}
			\xymatrix{ E_*\otimes E_* \ar[r]^(.6){\varphi_*} \ar[d]_{\theta_*\otimes \theta_*} & L \\
				E'_*\otimes E'_* \ar[ur]_(.6){\varphi'_*}
			}
		\end{align}
	\end{definition}
	
	We now describe the notion of parabolic stability and parabolic semi-stability for a parabolic symplectic vector bundle $(E_*,\ \varphi_*)$. Assume that the underlying vector bundle $E$ is of rank $r$ and degree $d$. Define the \textit{parabolic slope} of $E_*$
	to be
	\begin{align}\label{eqn:parabolic-slope}
		\mu_{par}(E_*)\,\,\,:=\,\,\, \dfrac{d+\sum_{p\in D}\sum_{i=1}^{\ell(p)}m^i_{_{p}}
			\alpha^i_{_{p}}}{r}\,\, \in\,\, {\mathbb R}.
	\end{align}
	Any algebraic sub-bundle $F$ of the underlying vector bundle $E$ gets equipped with an induced parabolic
	structure by restricting the flags and weights of $E_*$ to $F$. Let $F_*$ denote the resulting parabolic bundle.
	
	\begin{definition}[{\text{see also}\ \cite[Definition 
			2.1]{BMW11}}]\label{def:parabolic-symplectic-semistable-bundle}\mbox{}
		
		\begin{enumerate}
			\item Let $(E_*,\,\varphi_*)$ be a parabolic symplectic vector bundle (see
			Definition \ref{def:parabolic-symplectic-bundle}). An algebraic sub-bundle $F$ of the underlying bundle $E$ is
			said to be \textit{isotropic} if $\varphi_0(F\otimes F) \,=\, 0$; here $\varphi_0$ is the restriction of $\varphi_*$
			to $E_0\otimes E_0$. 
			
			\item $(E_*,\,\varphi_*)$ is said to be \textit{semistable parabolic symplectic}
			(respectively, \textit{stable parabolic symplectic}) if for all isotropic
			sub-bundles $0\, \not=\, F\,\subset\, E$ we have
			\begin{align*}
				\mu_{par}(F_*)\,\leq\, (\text{respectively,}\ \ <)\ \ \mu_{par}(E_*),
			\end{align*}
			where $F_*$ has the induced parabolic structure mentioned earlier.
			
			\item $(E_*,\ \varphi_*)$ is said to be a \textit{regularly stable parabolic symplectic vector bundle} if it is a
			stable parabolic symplectic bundle with the property that any nonzero (meaning not identically zero) parabolic endomorphism of
			$(E_*,\ \varphi_*)$ (see Definition \ref{def:parabolic-symplectic-isomorphism}) is
			multiplication by $\pm 1$.
		\end{enumerate}
	\end{definition}
	
	The maximal parabolic subgroups of the symplectic group ${\rm Sp}(r, {\mathbb C})$ are precisely those that preserve
	an isotropic subspace of ${\mathbb C}^r$ for the standard action of ${\rm Sp}(r, {\mathbb C})$ on ${\mathbb C}^r$. For
	this reason only the isotropic sub-bundles are used in Definition \ref{def:parabolic-symplectic-semistable-bundle}.
	
	\subsection{Semistable symplectic vector bundles in the non-parabolic case}\label{subsection:non-parabolic-case}\hfill\\
	
	Some results on the moduli space of usual (non-parabolic) semistable symplectic vector bundles on a curve	will be needed in order to compute the Brauer group of the parabolic symplectic moduli space. Therefore, for convenience, the relevant definitions for the non-parabolic case are recalled.
	Even though, by Remark \ref{rem:special-structure}, the non-parabolic case can equivalently be thought
	of as being endowed with the special trivial parabolic structure, it helps to discuss them 
	separately to avoid notational confusions.
	
	As before, fix a smooth irreducible complex projective curve $X$ of genus $g$, with $g\,\geq\, 3$. Fix an
	even integer $r\,\geq\, 2$ and a line bundle $L$ on $X$. Take a vector bundle $E$ on $X$ of rank $r$. Let
	$$
	\varphi\ :\ E\otimes E\ \longrightarrow\ L
	$$
	be a skew-symmetric $\mc{O}_X$--linear morphism of vector bundles.
	Let 
	\begin{equation}\label{e2-2}
		\widehat{\varphi} \,\,:\,\, E\,\,\longrightarrow\,\, E^{\vee}\otimes L\,=\, L\otimes E^\vee
	\end{equation}
	be the morphism defined by the following composition of maps:
	\begin{align*}
		E \,\simeq\, E\otimes \mc{O}_X\,\hookrightarrow\, E\otimes (E\otimes E^{\vee}) \,=\,
		(E\otimes E)\otimes E^{\vee} \,\xrightarrow{\,\,\, \varphi\otimes \text{Id}\,\,}\, L\otimes E^{\vee}.
	\end{align*}
	
	\begin{definition}\label{def:symplectic-bundle}
		A \textit{symplectic vector bundle} on $X$ taking values in $L$ is a pair $(E,\,\varphi)$ as above such that
		$\widehat{\varphi}$ in \eqref{e2-2} is an isomorphism.
	\end{definition}
	The notions of semi-stability and stability for symplectic vector bundles are standard
	(see \cite{Ra}). 
	
	\section{Fixed-point locus of parabolic symplectic moduli}\label{section:fixed-point-locus}
	
	As before, fix a smooth irreducible complex projective curve $X$ of genus $g$, with $g\,\geq\, 3$. Fix an 
	even integer $r\,\geq\, 2$ and a line bundle $L$ on $X$. Let $\mc{M}^{\boldsymbol{m,\alpha}}_L$ denote the moduli space of semistable 
	parabolic symplectic vector bundles $(E_*,\, \varphi_*)$ on $X$ of rank $r$ having a system of parabolic 
	weights $\boldsymbol{\alpha}$ and multiplicities $\boldsymbol{m}$, such that the symplectic form 
	$\varphi_*$ takes values in the line bundle $L$. The group of $2$--torsion line bundles on $X$ act on 
	$\mc{M}^{\boldsymbol{m,\alpha}}_L$ by tensor product. To describe this action, take a nontrivial line 
	bundle $\eta$ on $X$ of order two, and fix an isomorphism
	\begin{equation}\label{rho} \rho \,:\, \eta^{\otimes 2}\, \stackrel{\simeq}{\longrightarrow}\, \mc{O}_X.
	\end{equation}
	The action of
	$\eta$ sends a parabolic symplectic vector bundle $(E_* ,\, \varphi_*)$ to $(E_*\otimes \eta ,\,
	\varphi_*\otimes \rho)$. Thus, a
	fixed point under the action of $\eta$ is a parabolic symplectic vector bundle $(E_* ,\, \varphi_*)$
	admitting an isomorphism of parabolic symplectic vector bundles on $X$
	$$\theta_* \ :\ ( E_* , \, \varphi_*)\ \stackrel{\simeq}{\longrightarrow}\
	(E_*\otimes \eta ,\, \varphi_*\otimes\rho)$$
	such that the diagram
	\begin{align}\label{diagram:fixed-point}
		\xymatrix{E_*\otimes E_* \ar[rr]^{\varphi_*} \ar[d]_{\theta_*\otimes \theta_*} && L \\
			E_*\otimes E_*\otimes \eta^2\ar[urr]_{\varphi_*\otimes \rho} &&
		}
	\end{align}
	is commutative (see Definition \ref{def:parabolic-symplectic-isomorphism}). 
	
	We will describe the fixed-point locus of the semistable moduli
	space $\mc{M}^{\boldsymbol{m,\alpha}}_L$ for the action of a
	nontrivial $2$--torsion line bundle $\eta$. Using the isomorphism $\rho\,:\, \eta^{\otimes 2}\, \stackrel{\simeq}{\longrightarrow}\, \mc{O}_X$, the constant
	function $1$ on $X$ gives a nowhere vanishing section $s_0$ of $\eta^{\otimes2}$. Consider
	\begin{equation}\label{ey}
		Y \ :=\ \{y\, \in\, \eta\,\,\big\vert\,\, y^{\otimes 2}\, \in\, s_0(X)\}.
	\end{equation}
	The variety $Y$ thus constructed is an irreducible smooth projective curve. The restriction
	\begin{equation}\label{ega}
		\gamma \ :\ Y\ \longrightarrow\ X
	\end{equation}
	of the natural projection $\eta \, \longrightarrow\, X$ is actually a nontrivial \'etale Galois covering of
	degree two. There are well-defined notions of parabolic push--forward and pull--back of a parabolic vector
	bundle under a finite morphism; see \cite{AB23,BM19} for the details.
	
	Let $F_*$ be a parabolic vector bundle on $Y$ together with a skew-symmetric form
	\begin{equation}\label{eps}
		\phi_*\ :\ F_*\otimes F_*\ \longrightarrow \ \gamma^*L
	\end{equation}
	(see \eqref{ega} for $\gamma$). Consider the direct image of $\phi_*$
	\begin{equation}\label{ae}
		\gamma_*(\phi_*)\ :\ \gamma_*(F_*\otimes F_*)\ \longrightarrow\ \gamma_*\gamma^*L.
	\end{equation}
	The projection formula gives a homomorphism
	\begin{equation}\label{etr}
		\gamma_*\gamma^*L\ =\ L\otimes \gamma_*{\mathcal O}_Y \
		\xrightarrow{\,\,\,{\rm Id}_L\otimes\nu\,\,\, }\ L\otimes {\mathcal O}_X\ =\ L,
	\end{equation}
	where $\nu\, :\, \gamma_*{\mathcal O}_Y \,\longrightarrow\,{\mathcal O}_X$ is the trace map.
	Let $\widetilde{\nu}\, :\, \gamma_*\gamma^*L\,\longrightarrow \, L$ be the composition of maps in \eqref{etr}.
	Post--composing $\widetilde{\nu}$ with $\gamma_*(\phi_*)$ in \eqref{ae}, we have the homomorphism
	Using the homomorphism
	\begin{equation}\label{e3}
		\widetilde{\nu}\circ \gamma_*(\phi_*)\ :\ \gamma_*(F_*\otimes F_*)\ \longrightarrow\ L .
	\end{equation}
	Note that
	\begin{equation}\label{e4}
		\gamma^*(\gamma_* F_*)\ \simeq\ F_* \oplus \sigma^* F_*
	\end{equation}
	where $\sigma$
	is the (unique) nontrivial deck transformation of the degree two Galois covering $\gamma$. Now the
	projection map
	$$
	\gamma^*((\gamma_* F_*)^{\otimes 2})\,=\, (\gamma^*\gamma_* F_*)\otimes (\gamma^*\gamma_* F_*)\,
	\simeq\, (F_* \oplus \sigma^* F_*)\otimes (F_* \oplus \sigma^* F_*)
	\, \longrightarrow\, F_*\otimes F_*
	$$
	produces, using adjunction, a surjective homomorphism
	\begin{align}\label{eqn:adjunction-morphism}
		\iota_{_{F_*}} \ :\ (\gamma_* F_*)\otimes (\gamma_* F_*)\ \longrightarrow\ \gamma_* (F_*\otimes F_*);
	\end{align}
	see \cite[p.~110]{Ha} for adjunction. Combining this with the homomorphism in \eqref{e3}, we get a homomorphism
	\begin{align}\label{eqn:eq-1}
		\phi'_{*}\ :=\ (\widetilde{\nu}\circ \gamma_*(\phi_*))\circ \iota_{_{F_*}} \ :\ 
		\gamma_*(F_*)\otimes \gamma_*(F_*)\ \longrightarrow\ L .
	\end{align}
	It is evident from its construction that $\phi'_{*}$ is a skew-symmetric bilinear form on $\gamma_*(F_*)$ (recall
	that $\phi_*$ in \eqref{eps} is skew-symmetric).
	
	\begin{lemma}\label{lem:non-degenerate}
		For any semistable parabolic symplectic vector bundle $(F_*,\ \phi_*)$ on $Y$ taking values in the
		line bundle $\gamma^*L$, the pair $(\gamma_*(F_*),\ \phi'_*)$ in \eqref{eqn:eq-1}
		is a semistable parabolic symplectic vector bundle on $X$ taking values in $L$.
	\end{lemma}
	
	\begin{proof}
		It will be shown that there is a natural isomorphism of parabolic vector bundles
		\begin{align}\label{eqn:dual-commutes-with-pushroward}
			\Psi_*\ :\ \gamma_*(F_*^{\vee})\ \stackrel{\simeq}{\longrightarrow}\ (\gamma_* F_*)^{\vee} .
		\end{align}
		To see this, observe that $\gamma^*\left(\gamma_*(F_*^{\vee})\right)\,\simeq\, F_*^{\vee} \oplus \sigma^*(F_*^{\vee})\,\simeq\,
		F_*^{\vee}\oplus \sigma^*(F_*)^{\vee}$, while
		$$
		\gamma^*\left((\gamma_*F_*)^{\vee}\right) \,\simeq\, \left(\gamma^*\gamma_*(F_*)\right)^{\vee} \,\simeq\,
		\left(F_*\oplus\sigma^*(F_*)\right)^{\vee}\, \simeq\, F_*^{\vee}\oplus \sigma^*(F_*)^{\vee}.
		$$
		Thus, there is a natural parabolic isomorphism
		\begin{equation}\label{ei2}
			\gamma^*\left(\gamma_*(F_*^{\vee})\right)\ \stackrel{\simeq}{\longrightarrow}\ 
			\gamma^*\left((\gamma_*F_*)^{\vee}\right).
		\end{equation}
		It is straightforward to check that the isomorphism in \eqref{ei2} is actually equivariant for the actions of
		the Galois group ${\rm Gal}(\gamma)\,=\, {\mathbb Z}/2{\mathbb Z}$ on both sides. Therefore,
		the isomorphism in \eqref{ei2} descends to a parabolic isomorphism $\Psi_*$ as in \eqref{eqn:dual-commutes-with-pushroward}.
		
		Next, the non-degeneracy of $\phi_*$ implies the existence of the following parabolic isomorphism (Definition \ref{def:parabolic-symplectic-isomorphism}):
		$$
		\widehat{\phi}_*\ :\ F_*\ \stackrel{\simeq}{\longrightarrow}\ F_*^{\vee}\otimes \gamma^* L.
		$$
		Apply $\gamma_*$ to this isomorphism. Using projection formula and \eqref{eqn:dual-commutes-with-pushroward},
		we obtain a parabolic isomorphism
		\begin{align*}
			\gamma_*(F_*) \ \stackrel{\simeq}{\longrightarrow}\ (\gamma_* F_*)^{\vee}\otimes L, 
		\end{align*}
		which can be easily seen to coincide with the parabolic morphism $\widehat{\phi'_*}$ associated to
		$\phi'_{*}$ in \eqref{eqn:eq-1} (see \eqref{e2} for the construction of $\widehat{\phi'_*}$). It thus
		follows that $\widehat{\phi'_*}$ is an isomorphism, and hence $\phi'_*$ is non-degenerate. In other
		words, $(\gamma_*(F_*) ,\, \phi'_*)$ obtained using \eqref{eqn:eq-1} is a parabolic symplectic
		vector bundle taking values in the line bundle $L$.
		
		Assume that $(F_*,\ \phi_*)$ is a semistable parabolic symplectic vector bundle.
		To show the semistability of $(\gamma_*(F_*),\ \phi'_*)$, note that the condition that $(F_*,\ \phi_*)$ is semistable parabolic symplectic implies
		that $F_*$ is parabolic semistable \cite[Proposition 5.6]{BMW11}. It follows that $\gamma_*(F_*)$ is also parabolic semistable \cite[Proposition 4.3]{AB23}, and
		thus $\left(\gamma_*(F_*),\ \phi'_*\right)$ is a semistable parabolic symplectic vector bundle, again by \cite[Proposition 5.6]{BMW11}.
	\end{proof}
	
	It will be shown that $(\gamma_*(F_*),\, \phi'_*)$ is a fixed point for the action of $\eta$
	on the moduli space $\mc{M}_L^{\boldsymbol{m,\alpha}}$ that was described earlier.
	
	\begin{lemma}\label{lem:fixed-point}
		Fix an isomorphism $\rho\,:\, \eta^{\otimes 2}\,\stackrel{\simeq}{\longrightarrow}\,
		\mc{O}_X$ as in \eqref{rho}. The semistable symplectic parabolic vector bundle
		$(\gamma_*(F_*), \, \phi'_*)$ constructed in Lemma \ref{lem:non-degenerate} is a fixed point
		for the action of the line bundle $\eta$ on the semistable moduli space
		$\mc{M}_L^{\boldsymbol{m,\alpha}}$.
	\end{lemma}
	
	\begin{proof}
		Let $\pi\,:\, \eta\, \longrightarrow\, X$ be the natural projection. As $Y$ avoids the zeros of the tautological section of
		$\pi^*\eta$, the restriction of this tautological section produces a natural isomorphism 
		\begin{equation}\label{et}
			\tau\ :\ \mc{O}_Y\ \stackrel{\simeq}{\longrightarrow}\ \gamma^*\eta.
		\end{equation}
		Consider the isomorphism $\Id_{_{F_*}}\otimes \tau \,:\, F_*\, \stackrel{\simeq}{\longrightarrow}\, F_*\otimes \gamma^*\eta$.
		Applying $\gamma_*$ and using the projection formula $\gamma_*(F_*\otimes
		\gamma^*\eta)\,\simeq\, \gamma_*(F_*)\otimes \eta$ we obtain
		the following isomorphism of parabolic vector bundles on $X$:
		\begin{align}\label{eqn:fixed-point-map}
			\theta_*\ :\ \gamma_*F_* \ \stackrel{\simeq}{\longrightarrow}\ (\gamma_*F_*)\otimes \eta .
		\end{align}
		
		We will show that $\theta_*$ makes the diagram in \eqref{diagram:fixed-point} commutative
		for $(\gamma_*(F_*) ,\, \phi'_*)$, which would prove the lemma.
		
		First note that from the construction of the spectral curve $Y$ it is evident that
		\begin{align}\label{eqn:eq-2}
			\gamma^*\rho\ =\ (\tau^{\otimes 2})^{-1} \ :\ \gamma^*(\eta^{\otimes 2})\ 
			\stackrel{\simeq}{\longrightarrow}\ \mc{O}_Y,
		\end{align}
		where $\rho$ is the isomorphism in the lemma. This immediately yields the commutative diagram
		\begin{align}\label{diagram:diag-2}
			\xymatrix{
				F_*\otimes F_* \ar[rr]^{\phi_*} \ar[d]_{\Id_{_{F_*\otimes F_*}}\otimes \tau^2} && \gamma^*L\\
				F_*\otimes F_*\otimes \gamma^*\eta^2 \ar[urr]_{ \phi_*\otimes \gamma^*\rho}
			}
		\end{align}
		Apply $\gamma_*$ to this diagram. This produces a part of the following bigger commutative diagram:
		\begin{align*}
			\xymatrix@=2.1em{
				\gamma_*(F_*)\otimes \gamma_*(F_*) \ar@/^2.0pc/[rrrrrr]^{\phi'_*} \ar@{->>}[rr]^{\iota_{_{F_*}}}_{\eqref{eqn:adjunction-morphism}} \ar[dd]^{\simeq}_{\gamma_*(\Id_{_{F_*}}\otimes \tau)\otimes \gamma_*(\Id_{_{F_*}}\otimes \tau)} && \gamma_*(F_*\otimes F_*) \ar[rrr]^{\gamma_*(\phi_*)} \ar[dd]^{\simeq}_{\gamma_*((\Id_{_{F_*\otimes F_*}})\otimes \tau^2)} &&& \gamma_*\gamma^*L \ar[r] & L \\ &&&&&&\\
				\gamma_*(F_*\otimes\gamma^*\eta)\otimes\gamma_*(F_*\otimes\gamma^*\eta) \ar[rr]_(.55){\iota_{_{(F_*\otimes \gamma^*\eta)}}} \ar[dd]^{\simeq}_{\textnormal{proj. formula}} && \gamma_*(F_*\otimes F_*\otimes \gamma^*(\eta^{\otimes 2})) \ar[uurrr]^{\gamma_*(\phi_*\otimes\gamma^*\rho)} \ar[dd]^{\simeq}_{\textnormal{proj. formula}} &&&&\\&&&&&&\\
				(\gamma_*(F_*)\otimes \eta)\otimes(\gamma_*(F_*)\otimes \eta) \ar[rr]^(.55){\iota_{_{F_*}}\otimes \Id_{\eta^2}} && \gamma_*(F_*\otimes F_*)\otimes \eta^{\otimes 2} \ar[uuuurrr]_{\gamma_*(\phi_*)\otimes \rho}
			}
		\end{align*}
		Note that the composition of the two left-most vertical arrows in the above diagram is
		precisely $\theta_*\otimes \theta_*$ by construction. Thus, the outer-most arrows in the big diagram above give
		rise to the following diagram:
		\begin{align*}
			\xymatrix{
				\gamma_*(F_*)\otimes \gamma_*(F_*) \ar[rr]^(.6){\phi'_*} \ar[d]_{\theta_*\otimes \theta_*} && L\\
				\gamma_*(F_*)\otimes \gamma_*(F_*)\otimes \eta^2 \ar[urr]_(.6){\phi'_*\otimes \rho}
			}
		\end{align*}
		and hence $\theta_*$ makes the diagram in \eqref{diagram:fixed-point} commutative for
		$(\gamma_*(F_*) ,\, \phi'_*)$.
	\end{proof}
	
	Next, we discuss another construction through which parabolic symplectic vector bundles on $X$ can arise from
	the spectral curve $Y$, which shall also be fixed points for the action of $\eta$ on the moduli
	space $\mc{M}_L^{\boldsymbol{m,\alpha}}$. Let $F_*$ be a parabolic vector bundle on $Y$. Recall the decomposition in \eqref{e4}. Consider
	the following composition of homomorphisms constructed using \eqref{e4}:
	$$
	\gamma^*((\gamma_* F_*)^{\otimes 2})\,=\, (\gamma^*\gamma_* F_*)\otimes (\gamma^*\gamma_* F_*)\,
	=\, (F_* \oplus \sigma^* F_*)\otimes (F_* \oplus \sigma^* F_*)
	\, \longrightarrow\, F_*\otimes \sigma^* F_*.
	$$
	Using adjunction it produces a homomorphism
	\begin{equation}\label{e5}
		\mathcal{J}_{_{F_*}}\ : \ (\gamma_* F_*)\otimes (\gamma_* F_*)\ \longrightarrow \gamma_* (F_* \otimes \sigma^* F_*).
	\end{equation}
	
	Now let $$\psi_* \ :\ F_*\otimes \sigma^* F_* \ \longrightarrow\ \gamma^*L$$ be a
	pairing which is non-degenerate in the sense that the homomorphism
	\begin{equation}\label{ei}
		F_*\, \longrightarrow\, (\gamma^*L)\otimes (\sigma^* F_*)^\vee\,=\, (\gamma^*L)\otimes (\sigma^* F^\vee_*)
	\end{equation}
	induced by $\psi_* $ is an isomorphism of parabolic vector bundles. To describe another condition on
	$\psi_*$ that will be imposed, consider
	$$
	\sigma^*\psi_*\, :\, \sigma^*(F_*\otimes \sigma^* F_*)\,=\, (\sigma^* F_*)\otimes F_* \, \longrightarrow\,
	\sigma^*\gamma^*L\,=\, \gamma^*L.
	$$
	Assume that $\psi_*$ satisfies the following condition:
	\begin{align}\label{eqn:extra-condition}
		\sigma^*\psi_*\ =\ -\psi_*\circ \varpi,
	\end{align}
	where $\varpi\, :\, (\sigma^* F_*)\otimes F_*\, \longrightarrow\, F_*\otimes \sigma^* F_*$ is the natural involution that switches
	the two factors of the tensor product.
	
	As before, using the homomorphism $\mathcal{J}_{_{F_*}}$ in \eqref{e5} and the homomorphism
	$\gamma_*\gamma^* L\,\longrightarrow\, L$ in \eqref{etr}, consider the following composition of homomorphisms:
	$$
	(\gamma_* F_*)\otimes (\gamma_* F_*)\ \xrightarrow{\,\,\, \mathcal{J}_{_{F_*}}\,\,\,} \ \gamma_* (F_* \otimes \sigma^* F_*)
	\ \xrightarrow{\gamma_*(\psi_*)}\ \gamma_*\gamma^* L\ \longrightarrow\ L .
	$$
	We denote the composition of these morphisms by
	\begin{equation}\label{e6}
		\psi'_*\ :\ (\gamma_* F_*)\otimes (\gamma_* F_*)\ \longrightarrow \ L .
	\end{equation}
	
	\begin{lemma}\label{lem:another-fixed-point}
		Let $F_*$ be a parabolic semistable bundle on $Y$ together with a non-degenerate pairing
		$$\psi_*\ :\ F_*\otimes (\sigma^*F_*)\ \longrightarrow\ \gamma^*L$$ satisfying 
		\eqref{eqn:extra-condition}. Then the following two statements hold: 
		\begin{enumerate}[(1)]
			\item $(\gamma_*F_*, \ \psi'_*)$ \textnormal{(}$\psi'_*$ as defined in \eqref{e6}\textnormal{)} is a
			semistable parabolic symplectic vector bundle on $X$ taking values in the line bundle $L$.
			
			\item Fix an isomorphism $\rho\,:\, \eta^{\otimes 2}\,\stackrel{\simeq}{\longrightarrow}\, \mc{O}_X$
			as in \eqref{rho}. The semistable parabolic symplectic vector bundle $(\gamma_*(F_*), \, \psi'_*)$ in (1) is a
			fixed point for the action of the line bundle $\eta$ on the semistable moduli space $\mc{M}^{\boldsymbol{m,\alpha}}_L$.
		\end{enumerate}
	\end{lemma}
	
	\begin{proof}
		The proof is along the same lines as that of Lemma \ref{lem:fixed-point} with some modifications; for the sake of completeness we provide the details.
		
		Proof of (1)\ : \ Since $\sigma$ is an isomorphism, both $\sigma^*\circ\sigma_*$ and $\sigma_*\circ\sigma^*$
		induce the identity functor. This, combined with the fact that $\sigma$ is an involution, gives
		$$\sigma^*\ =\ (\sigma_*)^{-1}\ =\ \sigma_*.$$
		Next, the non-degeneracy of $\psi_*$ gives rise to a parabolic isomorphism
		\begin{equation}\label{ez}
			\widehat{\psi_*} \ :\ F_*\ \stackrel{\simeq}{\longrightarrow} \
			(\sigma^*F_*)^{\vee}\otimes \gamma^*L\ =\ (\sigma_*F_*)^{\vee}\otimes \gamma^*L 
		\end{equation}
		(see \eqref{ei}). Applying $\gamma_*$ on both sides of \eqref{ez} and using the projection formula,
		$$
		\gamma_*(\widehat{\psi_*}) \ :\ \gamma_*(F_*)\ \longrightarrow \ \gamma_*((\sigma_*F_*)^{\vee})\otimes L.
		$$
		Now using the isomorphism in \eqref{eqn:dual-commutes-with-pushroward} and the equality
		$\gamma\circ\sigma\,=\,\gamma$, we obtain an isomorphism of parabolic vector bundles
		\begin{equation}\label{ez2}
			\gamma_*(F_*) \ \stackrel{\simeq}{\longrightarrow}\ \gamma_*(F_*)^{\vee}\otimes L.
		\end{equation}
		It is straightforward to check that the isomorphism in \eqref{ez2}
		coincides with the parabolic morphism $\widehat{\psi'_*}$ constructed using $\psi'_*$ just as
		it was done in \eqref{e2}. Consequently, $\psi'_*$ is non-degenerate. Moreover, the condition
		imposed on $\psi_*$ in \eqref{eqn:extra-condition} ensures that $\psi'_*$
		is skew-symmetric as well.
		
		The fact that $(\gamma_*F_*, \ \psi'_*)$ is semistable follows from the same argument as the one given in the proof of	Lemma \ref{lem:non-degenerate}.
		
		Proof of (2)\ : \ Consider the isomorphism $\theta_*\,:\, \gamma_*(F_*)\,\stackrel{\simeq}{\longrightarrow}
		\,\gamma_*(F_*)\otimes \eta$ in \eqref{eqn:fixed-point-map}. Now, using \eqref{eqn:eq-2} we have the following
		diagram which is analogous to \eqref{diagram:diag-2}:
		\begin{align*}
			\xymatrix{
				F_*\otimes (\sigma^*F_*) \ar[rr]^{\psi_*} \ar[d]_{\Id_{_{F_*\otimes (\sigma^*F_*)}}\otimes \tau^2} && \gamma^*L\\
				F_*\otimes (\sigma^*F_*)\otimes \gamma^*\eta^2 \ar[urr]_{ \psi_*\otimes \gamma^*\rho}
			}
		\end{align*}
		Apply $\gamma_*$ to this diagram produces a part of the following bigger diagram similar to
		the one in Lemma \ref{lem:fixed-point} (note that we have to use $\gamma\circ\sigma
		\,=\, \gamma$ at some places in the diagram below):
		\begin{align*}
			\xymatrix@=2.1em{
				\gamma_*(F_*)\otimes \gamma_*(F_*) \ar@/^2.0pc/[rrrrrr]^{\psi'_*} \ar[rr]_{\eqref{e5}}^{\mathcal{J}_{_{F_*}}} \ar[dd]_{\gamma_*(\Id_{_{F_*}}\otimes \tau)\otimes \gamma_*(\Id_{_{F_*}}\otimes \tau)} && \gamma_*(F_*\otimes \sigma^*F_*) \ar[rrr]^{\gamma_*(\psi_*)} \ar[dd]_{\gamma_*((\Id_{_{F_*\otimes (\sigma^*F_*)}})\otimes \tau^2)} &&& \gamma_*\gamma^*L \ar[r] & L \\ &&&&&&\\
				\gamma_*(F_*\otimes\gamma^*\eta)\otimes\gamma_*(F_*\otimes\gamma^*\eta) \ar[rr]_(.55){\mathcal{J}_{_{(F_*\otimes \gamma^*\eta)}}} \ar[dd]^{\simeq}_{\textnormal{proj. formula}} && \gamma_*(F_*\otimes (\sigma^*F_*)\otimes \gamma^*(\eta^{\otimes 2})) \ar[uurrr]^{\gamma_*(\psi_*\otimes\gamma^*\rho)} \ar[dd]^{\simeq}_{\textnormal{proj. formula}} &&&&\\&&&&&&\\
				(\gamma_*(F_*)\otimes \eta)\otimes(\gamma_*(F_*)\otimes \eta) \ar[rr]^(.55){\mathcal{J}_{_{F_*}}\otimes \Id_{\eta^2}} && \gamma_*(F_*\otimes \sigma^*F_*)\otimes \eta^{\otimes 2} \ar[uuuurrr]_{\gamma_*(\psi_*)\otimes \rho}
			}
		\end{align*}
		The outer-most arrows of this big diagram produces the following diagram:
		\begin{align*}
			\xymatrix{
				\gamma_*(F_*)\otimes \gamma_*(F_*) \ar[rr]^(.6){\psi'_*} \ar[d]_{\theta_*\otimes \theta_*} && L\\
				\gamma_*(F_*)\otimes \gamma_*(F_*)\otimes \eta^2 \ar[urr]_(.6){\psi'_*\otimes \rho}
			}
		\end{align*}
		which implies that $(\gamma_*F_*, \ \psi'_*)$ is a fixed point for the action of $\eta$
		on the moduli space $\mc{M}^{\boldsymbol{m,\alpha}}_L$.
	\end{proof}
	\vspace{0.2cm}
	
	\subsection{Codimension of the fixed point locus for full-flag systems}
	\hfill\\
	
	In this subsection it will be assumed that the parabolic structure consists of \textit{full-flag}
	systems at each of the parabolic points, meaning $m^i_{_{x}} \,=\,1$ for all $x\,\in\, S$ and $i$
	(see Definition \ref{def:parabolic-bundles}). Consequently, we shall drop the
	symbol $'\boldsymbol{m}'$ from the notation of the parabolic symplectic moduli in this section. More precisely, fix an even positive integer $r$
	and a line bundle $L$ on $X$. Denote by $\mc{M}^{\boldsymbol{\alpha}}_L$ the moduli space of semistable parabolic symplectic vector
	bundles $(E_*,\, \varphi_*)$ on $X$ of rank $r$ and having a full-flag system of multiplicities
	and a system of parabolic weights $\boldsymbol{\alpha}$; the symplectic form
	$\varphi_*$ takes values in the line bundle $L$. Since $(E_*,\, \varphi_*)$ is a parabolic
	symplectic bundle, the systems of weights $\boldsymbol{\alpha}$
	are of symmetric type, in the sense of \cite[Definition 3.4]{BCD24}. This symmetry condition
	is described below.
	
	Let
	\begin{equation*}
		S\ \subset\ X
	\end{equation*}
	be the set of parabolic points. For each point $x\, \in\, S$, fix real numbers
	$$
	0\, <\, \alpha^1_{_{x}} \, <\, \alpha^2_{_{x}} \, <\, \cdots \, <\, \alpha^{{{r/2}}}_{_{x}}\, < \, \frac{1}{2}.
	$$
	So the set of $r$--numbers $\alpha^1_{_{x}}\ ,\, \cdots , \, \alpha^{{r/2}}_{_{x}}\ ,\ 
	1-\alpha^1_{_{x}}\ ,\, \cdots,\, 1- \alpha^{{r/2}}_{_{x}}$ are all distinct. The systems of parabolic
	weights $\boldsymbol{\alpha}$ will be given by the sequence $$\left(\alpha^1_{_{x}}\ ,\ \alpha^2_{_{x}}\ ,\ \cdots\ ,\ \alpha^{r/2}_{_{x}}\ ,\ 1-\alpha^{r/2}_{_{x}}\ ,\ \cdots\ ,\ 1-\alpha^1_{_{x}}\right)$$ at each $x\in S$. Let 
	\begin{equation}\label{ers}
		\mc{M}^{{\boldsymbol{\alpha}},rs}_L\ \,\subset\ \,\mc{M}^{\boldsymbol{\alpha}}_L
	\end{equation}
	denote the Zariski open subset consisting of regularly stable parabolic symplectic vector
	bundles (see Definition \ref{def:parabolic-symplectic-semistable-bundle}); recall the convention
	that the symbol $'\boldsymbol{m}'$ is dropped. It is 
	straightforward to see that $\mc{M}^{\boldsymbol{\alpha},rs}_L$ is invariant under the 
	action, on $\mc{M}^{\boldsymbol{\alpha}}_L$, of the $2$--torsion line bundles over $X$.
	
	Given a nontrivial $2$--torsion line bundle $\eta$ on $X$, consider the spectral curve 
	$$\gamma\ :\ Y \ \longrightarrow \ X$$ constructed as in \eqref{ey}. We need to 
	describe certain moduli spaces of parabolic bundles on $Y$ with parabolic 
	structures along $\gamma^{-1}(S)$, so some conventions will be employed. For each 
	$x\,\in\, S$, set $$A_x \ :=\ \{\alpha^{1}_{_{x}}\ ,\,\alpha^{2}_{_{x}}\ 
	,\,\cdots,\,\alpha_{_{x}}^{r/2}\ ,\, 1-\alpha^{1}_{_{x}},\,\cdots,\, 
	1-\alpha^{r/2}_{_{x}}\}.$$ Let $P(A_x)$ denote the collection of all possible partitions 
	of $A_x$ into two disjoint subsets of cardinality $\frac{r}{2}$ each. In any such 
	partition, the two disjoint subsets shall be referred to as \textit{cells}. Let 
	$$P_1(A_x)\ \subset\ P(A_x)$$ be the subset consisting of all partitions with the 
	following property: If $\alpha_{_{x}}^{i}$ belongs to a cell, then 
	$1-\alpha^{i}_{_{x}}$ belongs to the \textit{same} cell. On the other hand, let 
	$$P_2(A_x)\ \subset\ P(A_x)$$ be the subset consisting of all partitions which 
	satisfy the property that if $\alpha^{i}_{_{x}}$ belongs to a cell, then 
	$1-\alpha^{i}_{_{x}}$ belongs to the \textit{other} cell.
	
	For simplicity, assume that we have a single parabolic point $S\,=\,\{x\}$. Consider $A_x$ 
	as above. Now, for each partition $\textbf{t}\,\in\, P_1(A_x)$, let 
	$\mc{M}^{\textbf{t}}_{Y}$ denote the moduli space of semistable parabolic symplectic vector 
	bundles of rank $\frac{r}{2}$ on $Y$ such that the symplectic form takes values in 
	$\gamma^*L$, with parabolic structure consisting of full flags along the points of 
	$\gamma^{-1}(x)$ and weights assigned using the partition $\textbf{t}$, as constructed in 
	detail in \cite[\S~3]{BCD23}. Note that $\mc{M}^{\textbf{t}}_{Y}$ is empty if $\frac{r}{2}$ 
	is an odd integer. Denote
	\begin{align*}
		\mc{M}_Y\ :=\ \coprod_{\textbf{t}\in P_1(A_x)}\mc{M}^{\textbf{t}}_{Y}.
	\end{align*} 
	Also, for each $\textbf{t}'\,\in\, P_2(A_x)$, let $\mc{N}^{\textbf{t}'}_Y$ denote the 
	moduli space of parabolic semistable vector bundles $F_*$ of rank $\frac{r}{2}$ on $Y$ 
	having full-flag parabolic structures along the points of $\gamma^{-1}(x)$ and parabolic 
	weights assigned using the partition $\textbf{t}'$, together with a pairing $\psi_*\, :\, 
	F_*\otimes\sigma^*F_* \,\longrightarrow \,\gamma^*L$ satisfying the condition in 
	\eqref{eqn:extra-condition}. Denote
	\begin{align*}
		\mc{M}'_Y\ :=\ \coprod_{\textbf{t}'\in P_2(A_x)}\mc{N}^{\textbf{t}'}_{Y}.
	\end{align*} 
	See the proof of \cite[Proposition 3.3]{BCD24} for the reason why only such partitions
	are allowed. Finally, when the cardinality of the set of parabolic points
	$S\, \subset\, X$ is greater than 1, we perform the same procedure for each parabolic point
	to obtain our system of parabolic weights on $\gamma^{-1}(S)$.
	
	Let $\left(\mc{M}^{\boldsymbol{\alpha}}_L\right)^{\eta} \ \subset\ \mc{M}^{\boldsymbol{\alpha}}_L$ as well as $\left(\mc{M}^{\boldsymbol{\alpha},rs}_L\right)^{\eta}\ \subset\ \mc{M}^{\boldsymbol{\alpha},rs}_L$ (see \eqref{ers}) denote the fixed-point loci for the action of $\eta$.
	By Lemma \ref{lem:fixed-point} and Lemma \ref{lem:another-fixed-point}, there exists a
	natural morphism \begin{align}\label{eqn:morphism-to-fixed-point}
		f\ :\ \mc{M}_{Y}\coprod \mc{M}'_Y\ \longrightarrow\ 
		\left(\mc{M}^{\boldsymbol{\alpha}}_L\right)^{\eta}.
	\end{align}
	
	\begin{proposition}\label{prop:fixed-point-morphism-surjective}
		Let $V \,:=\, f^{-1}\left((\mc{M}_L^{\boldsymbol{\alpha},rs})^{\eta}\right)$
		(see \eqref{ers}) with $f$
		being the map in \eqref{eqn:morphism-to-fixed-point}. The restricted morphism
		$$f\big\vert_V\ :\ V\ \longrightarrow\ (\mc{M}_{L}^{\boldsymbol{\alpha},rs})^{\eta}$$ is surjective.
	\end{proposition}
	
	\begin{proof}
		To show that $f\big\vert_V$ is surjective, take any $(E_*,\, \varphi_*)\,\in\, (\mc{M}_L^{\boldsymbol{\alpha},rs})^{\eta}$, so
		that there exists an isomorphism $\rho\,:\, \eta^2\,\xrightarrow{\,\,\,\simeq\,\,\,} \mc{O}_X$
		(as in \eqref{rho}) and an isomorphism of parabolic symplectic vector bundles
		$$\theta_* \ :\ (E_*,\, \varphi_*)\ \xrightarrow{\,\,\,\simeq\,\,\,}\
		(E_*\otimes \eta,\, \varphi_*\otimes \rho)$$ 
		in the sense of Definition \ref{def:parabolic-symplectic-isomorphism}, yielding the
		commutative diagram \eqref{diagram:fixed-point}. Observe that 
		\begin{align}\label{eqn:iso-1}
			\Id_{_{E_*}}\otimes \rho^{-1} \ :\ (E_*,\, \varphi_*)\ \xrightarrow{\,\,\,\simeq\,\,\,}\
			(E_*\otimes \eta^2,\, \varphi_*\otimes \rho^2)
		\end{align}
		is an isomorphism of parabolic symplectic vector bundles, where $\rho$ is the isomorphism
		in \eqref{rho}. Now, as $(E_*, \, \varphi_*)$
		is a parabolic regularly stable symplectic vector bundle (see Definition
		\ref{def:parabolic-symplectic-semistable-bundle}), any two parabolic automorphisms of $(E_*,\, \varphi_*)$ differ by
		multiplication with an element of $\pm$ 1. 
		Thus, we can re-scale $\theta_*$, if necessary, to ensure
		that the parabolic morphism
		$\theta_*\circ\theta_*$ coincides with $\Id_{_{E_*}}\otimes \rho^{-1}$ in \eqref{eqn:iso-1}.
		Now, the proof in \cite[Lemma 3.3]{BCD23}
		produces a parabolic vector bundle $F_*$ of rank $\frac{r}{2}$ on $Y$ such that $\gamma_*( F_*)
		\,\simeq\, E_*$. As the construction of the parabolic structure along $\gamma^{-1}(S)$ was done only for full-flag systems of weights in \cite[Lemma 3.3]{BCD23}, we have restricted ourselves to the full-flag situation here. We briefly recall the construction, for the sake of convenience of the reader.
		
		Composing $\gamma^* \theta$ with $\Id_{_{\gamma^*E}}\otimes \tau^{-1}$, where $\tau$ is the isomorphism in \eqref{et}, one obtains an endomorphism of $\gamma^* E$, which is denoted by
		\begin{equation}\label{tp}
			\theta'\ \in\ H^0(Y,\, \text{End}(\gamma^* E))\ =\ H^0(Y,\, \gamma^*\text{End}(E)).
		\end{equation}
		{}From the construction
		of $\theta'$ it follows that $\theta'\circ\theta'
		\,=\, \Id_{_{\gamma^*E}}$, and thus $\gamma^* E$ decomposes into a direct sum of sub-eigen-bundles corresponding to
		the two eigenvalues
		$\pm \,1$. We take $F$ to be the eigenbundle corresponding to the eigenvalue $1$, and equip it with the
		induced parabolic structure on $\gamma^{-1}(S)\, \subset\, Y$ to obtain $F_*$. 
		
		We have a parabolic isomorphism
		$$ \gamma^*E_* \ \simeq\ F_*\oplus (\sigma^*F_*),$$
		where $F_*$ and $\sigma^*F_*$ are the parabolic sub-bundles of $\gamma^*E_*$ corresponding
		to the sub-eigen-bundles of $\gamma^*E$ for the eigenvalues $1$ and $-1$ respectively. It now follows that 
		\begin{align*}
			(\gamma^*E_*)\otimes (\gamma^*E_*)\ \simeq\ (F_*\otimes F_*)\bigoplus (F_*\otimes
			(\sigma^*F)_*)\bigoplus ((\sigma^*F)_*\otimes(\sigma^*F)_*)\bigoplus ((\sigma^*F)_*\otimes F_*) .
		\end{align*}
		Set $V_1 \,:=\, (F_*\otimes F_*)\bigoplus \sigma^*(F_*\otimes F_*)\,=\,
		(F_*\otimes F_*)\bigoplus (\sigma^* F_*)\otimes (\sigma^*F_*)$ and
		$$V_{-1}\,:=\, (F_*\otimes(\sigma^*F)_*)\bigoplus
		((\sigma^*F)_*\otimes F_*).$$ Clearly, $V_1$ and $V_{-1}$ are the sub-eigen-bundles of
		$(\gamma^*E_*)\otimes (\gamma^*E_*)$ for the eigenvalues $1$ and $-1$ respectively. Moreover, both $V_1$ and $V_{-1}$ are equivariant
		sub-bundles for the action of the Galois group ${\rm Gal}(\gamma)\,=\,
		{\mathbb Z}/2{\mathbb Z}$ on $(\gamma^*E)_*$. It follows that the equivariant parabolic
		morphism $\gamma^*\varphi_*$ is completely determined by the two parabolic
		morphisms of the following form
		\begin{align}\label{eqn:upstairs-bilinear-forms}
			\zeta_*\,:\,F_*\otimes F_*\,\longrightarrow\, \gamma^*L \ \ \text{and}\ \ \xi_*\,:\,
			F_*\otimes(\sigma^*F)_*\,\longrightarrow\, \gamma^*L .
		\end{align}
		By uniqueness of descent it follows that
		$$E_*\otimes E_*\ \stackrel{\simeq}{\longrightarrow}\ \gamma_*(F_*\otimes F_*)\
		\bigoplus\ \gamma_*\left(F_*\otimes (\sigma^*F)_*\right);$$
		under this isomorphism, $\varphi_*$ corresponds to the parabolic morphism
		\begin{align}\label{eqn:tensor-decomposition}
			\gamma_*(\zeta_*) + \gamma_*(\xi_*)\ :\ \gamma_*(F_*\otimes F_*)\ \bigoplus\
			\gamma_*\left(F_*\otimes (\sigma^*F)_*\right)\ \longrightarrow\ L,
		\end{align}
		where $\zeta_*$ and $\xi_*$ are the homomorphisms in \eqref{eqn:upstairs-bilinear-forms}. 
		
		Next, we pre-compose $\gamma_*(\zeta_*)$ (respectively, $\gamma_*(\xi_*)$) with the usual parabolic morphism
		$$\gamma_*(F_*)\otimes \gamma_*(F_*)\,\underset{\eqref{eqn:adjunction-morphism}}{\longrightarrow} \,\gamma_*(F_*\otimes F_*)\ \ \left(\text{respectively,}\ \ 
		\gamma_*(F_*)\otimes \gamma_*\left((\sigma^*F)_*\right)\,\longrightarrow\,
		\gamma_*\left(F_*\otimes (\sigma^*F)_*\right)\right).$$ 
		The resulting parabolic
		morphisms thus obtained coincide with $\zeta'_*$ and $\xi'_*$ respectively, where $\zeta'_*$ and
		$\xi'_*$ are constructed as in \eqref{eqn:eq-1} and \eqref{e6}. Also, $F_*$ is parabolic
		semistable, because $E_*\, =\,\gamma_*(F_*)$ is parabolic semistable.
		
		To analyze \eqref{eqn:tensor-decomposition} further, note that
		due to the effect of re-scaling $\theta_*$ in the beginning of the proof, the diagram \eqref{diagram:fixed-point} may
		not commute anymore.
		Replace $\rho$ by $\lambda\cdot \rho$, where $\lambda$ is an appropriate scalar, so
		that the following diagram is commutative:
		\begin{align*}
			\xymatrix{E_*\otimes E_* \ar[rr]^{\varphi_*} \ar[d]_{\theta_*\otimes \theta_*} && L \\
				E_*\otimes E_*\otimes \eta^2\ar[urr]_{\varphi_*\otimes \lambda\rho} &&
			}
		\end{align*} 
		Apply $\gamma^*$ to it, and compose the left vertical arrow with $\Id_{_{(\gamma^*E_*\otimes\gamma^*E_*)}}\otimes\tau^{-2}$ to get the following:
		\begin{align}\label{diagram:diag-1}
			\xymatrix{(\gamma^*E)_*\otimes(\gamma^*E)_* \ar[rrr]^{\gamma^*(\varphi_*)} \ar[d]_{\gamma^*(\theta_*)\otimes \gamma^*(\theta_*)} &&& \gamma^*L \\
				(\gamma^*E)_*\otimes(\gamma^*E)_*\otimes\gamma^*\eta^2\ar[rrru]^{\gamma^*(\varphi_*)\otimes\lambda\gamma^*\rho} \ar[d]_{\Id_{_{(\gamma^*E_*\otimes\gamma^*E_*)}}\otimes\tau^{-2}}&&&\\
				(\gamma^*E)_*\otimes(\gamma^*E)_*\ar[rrruu]_{\gamma^*(\varphi_*)\otimes \lambda(\gamma^*\rho\otimes\tau^2)}&&&
			}
		\end{align}
		As mentioned earlier, $\gamma^*\rho\otimes\tau^2$ is the identity map on the trivial line
		bundle $\mc{O}_Y$. Also, note that the composition of the left vertical arrows
		in \eqref{diagram:diag-1} coincides with $\theta'_*\otimes\theta'_*$, where $\theta'_*\,:\, (\gamma^*E)_*\,
		\longrightarrow\, (\gamma^*E)_*$ is the composition of
		$\gamma^* \theta_*$ with $\Id_{_{\gamma^*E_*}}\otimes\tau^{-1}$ as
		in \eqref{tp}. Consequently, the diagram \eqref{diagram:diag-1} takes on the following form:
		\begin{align}\label{diagram:diag-3}
			\xymatrix{
				(\gamma^*E)_*\otimes(\gamma^*E)_* \ar[rr]^(.56){\gamma^*(\varphi_*)} \ar[d]_{\theta'_*\otimes \theta'_*} && \gamma^*L \\
				(\gamma^*E)_*\otimes(\gamma^*E)_*\ar[rru]_(.56){\lambda\cdot\gamma^*(\varphi_*)} &&
			}
		\end{align}
		
		To analyze the possible values of $\lambda$, consider the following two cases,
		depending on the behaviour of $F$ under the bilinear form $\gamma^*\varphi$.
		
		\textbf{Case I :}\, Suppose that $F$ is not an isotropic sub-bundle of $\gamma^*E$ under
		$\gamma^*\varphi$. Thus, there exists two nonzero vectors $v_1,\,v_2$ in $F$ with
		$\gamma^*(\varphi)(v_1,\,v_2)\,\neq\, 0$. As $F$ is the sub-eigen-bundle of $\gamma^*E$
		for the eigenvalue $1$ of the automorphism $\theta'$, it follows from the diagram
		\eqref{diagram:diag-3} that
		$$\gamma^*(\varphi)(v_1,\,v_2)\ =\ \lambda\cdot\gamma^*(\varphi)(\theta'(v_1),\,\theta'(v_2))
		\ =\ \lambda\cdot\gamma^*(\varphi)(v_1,\,v_2),$$
		and thus $\lambda\,=\,1$. Now, as $\sigma^* F_*\,\subset\, (\gamma^*E)_*$ is the sub-eigen-bundle 
		for the eigenvalue $-1$ of the parabolic automorphism $\theta'_*$ of $(\gamma^*E)_*$,
		it follows immediately from the
		diagram \eqref{diagram:diag-3} that $\sigma^*(F_*)$ is the orthogonal complement of
		$F_*$ under the form $\gamma^*(\varphi_*)$. As $\gamma^*(E_*)\,=\,
		F_* \oplus \sigma^*(F_*)$, it is deduced
		that the restriction of $\gamma^*(\varphi_*)$ to $F_*\otimes F_*$, namely $\zeta_*$, is non-degenerate on $F_*$, while the restriction of $\gamma^*(\varphi_*)$ to $F_*\otimes \sigma^*(F_*)$, namely $\xi_*$, is the
		zero map (see \eqref{eqn:upstairs-bilinear-forms}). In other words, $(F_*, \, \zeta_*)$
		is a parabolic semistable symplectic vector bundle, and it follows from
		\eqref{eqn:tensor-decomposition} that $(E_*, \, \varphi_*)\,\simeq\, \left(\gamma_*(F_*),\,
		\zeta'_{*}\right)$ as parabolic symplectic vector bundles, where $\zeta'_*$ is as in \eqref{eqn:eq-1}.
		
		\textbf{Case II :}\, $F$ is an isotropic sub-bundle of $\gamma^*(E)$ under $\gamma^*(\varphi)$. In this case, by the 
		non-degeneracy of $\gamma^*(\varphi_*)$, for any nonzero vector $v_1$ in $F$
		we can find a nonzero vector $v_2$ in $\sigma^*F$ satisfying the condition
		$\gamma^*(\varphi)(v_1,\,v_2)\,\neq\, 0$. Again, using diagram \eqref{diagram:diag-3} and
		the fact that $\sigma^*(F_*)$ is the sub-eigen-bundle of $\gamma^*(E_*)$, for eigenvalue
		$-1$ of $\theta'_*$, we get that
		$$\gamma^*(\varphi)(v_1,\,v_2)\ =\ \lambda\cdot\gamma^*(\varphi)\left(\theta'(v_1),\,
		\theta'(v_2)\right)\ =\ -\lambda\cdot\gamma^*(\varphi)(v_1,\,v_2),$$
		and thus $\lambda\,=\,-1$. In this case, the restriction of $\gamma^*(\varphi_*)$ to
		$F_*\otimes (\sigma^*F)_*$ is non-degenerate, and $\zeta_* \,=\,0$
		in \eqref{eqn:upstairs-bilinear-forms}. Moreover, from the fact that $\varphi_*$ is
		skew-symmetric, it follows that $\xi_*$ satisfies the equation \eqref{eqn:extra-condition}
		as well. Consequently, $(E_*,\, \varphi_*)\,\simeq\, (\gamma_*(F_*),\, \xi'_{*})$, where
		$\xi'_*$ is as in \eqref{e6}. This completes the proof of the proposition.
	\end{proof}
	
	\begin{corollary}\label{cor:free-action-complement-codimension}
		Let $\Gamma$ denote the finite group of $2$--torsion line bundles on $X$. Consider the
		Zariski closed subset 
		$$Z_{\boldsymbol{\alpha}}\ :=\ \bigcup_{\eta\in\Gamma\setminus\{\mc{O}_X\}}
		(\mc{M}_L^{\boldsymbol{\alpha},rs})^{\eta}\ \subset\ \mc{M}^{\boldsymbol{\alpha},rs}_L .$$
		The codimension of $Z_{\boldsymbol{\alpha}}$ in $\mc{M}^{\boldsymbol{\alpha},rs}_L$
		is at least $3$.
	\end{corollary}
	
	\begin{proof}
		Let $\eta\,\in\, \Gamma$ be a nontrivial line bundle. Let $\gamma\,:\,Y\,\longrightarrow\,
		X$ be the spectral curve corresponding to $\eta$ as before. To prove the result it suffices
		to show that
		\begin{equation}\label{es}
			\dim(\mc{M}^{\boldsymbol{\alpha},rs}_L)-
			\dim((\mc{M}^{\boldsymbol{\alpha},rs}_L)^{\eta})\ \geq\ 3.
		\end{equation}
		Now, by Proposition
		\ref{prop:fixed-point-morphism-surjective},
		$$\dim\left((\mc{M}_L^{\boldsymbol{\alpha},rs})^{\eta}\right)\ \leq\
		\max\left\{\dim(\mc{M}_Y),\ \dim(\mc{M}'_Y)\right\}.$$
		Let the cardinality of the set of parabolic points $S$ is $|S|\,=\,
		s$; so the cardinality of $\gamma^{-1}(S)$ is $|\gamma^{-1}(S)|\,=\,2s$. Let $g(Y)$ denote the
		genus of $Y$. By Riemann--Hurwitz formula, $g(Y) \,=
		\, 2(g-1)+1$. For notational convenience, denote $p\,:=\, \frac{r}{2}$. The dimension of full-flag symplectic isotropic flag varieties can be found using the formula given in \cite[\S~2]{Co14}, which turns out to be $p^2$. Then, using \cite[Lemma 3.10]{BR88}, we have the following expressions for	dimensions:
		\begin{align*}
			&\dim(\mc{M}^{\boldsymbol{\alpha},rs}_L) \,=\, \dim(\mc{M}^{\boldsymbol{\alpha}}_L)\,=\,
			p(2p+1)(g-1)+sp^2\ ,\ \ \left[s\,=\,|S|,\ p\,=\,\frac{r}{2}\right]\\
			\text{while}\ \ &\dim(\mc{M}_Y) \,=\, 0 \ \ \textnormal{(if}\ \ p \ \ \textnormal{is odd), and} \nonumber\\
			& \dim(\mc{M}_Y)\, =\,\dfrac{p}{2}(p+1) \left(2g(Y)-1\right)+2s\cdot
			\left(\frac{p}{2}\right)^2
			\ \ \textnormal{(if}\ p \ \textnormal{is even)}\nonumber\\
			& \qquad \qquad =\ \frac{p}{2}\left(p+1\right)(2g-2)+ 2s\cdot\frac{p^2}{4}\nonumber\\
			&\qquad\qquad =\ p(p+1)(g-1)+ s\cdot\frac{p^2}{2}.
		\end{align*}
		Also, notice that for a parabolic vector bundle $F_*$ in $\mc{M}'_Y$, the
		line bundle $\det(\gamma_*F)$ is fixed. Thus
		\begin{align*}
			\dim(\mc{M}'_Y) &\,\leq\, p^2\left(g(Y) -1\right) +1 -g + 2s\cdot \frac{p(p-1)}{2}\nonumber\\
			&=\, p^2(2g-2) +1 -g + s\cdot p(p-1)\nonumber\\
			&=\, (p^2-1)(g-1)+ s\cdot p(p-1).
		\end{align*}
		It follows that 
		\begin{align*}
			\dim(\mc{M}^{\boldsymbol{\alpha},rs}_L) - \dim(\mc{M}_Y)\, &
			=\, \left(p(2p+1)(g-1)+sp^2\right) - \left(p(p+1)(g-1)+ s\cdot\frac{p^2}{2}\right)\nonumber\\
			& =\,p^2(g-1) + s\cdot \dfrac{p^2}{2}\nonumber\\
			&\geq\ 3 .
		\end{align*}
		\begin{align*}
			\text{Similarly},\, \ \dim(\mc{M}^{\boldsymbol{\alpha},rs}_L)-
			\dim(\mc{M}'_Y) \,&\geq\, \left(p(2p+1)(g-1)+sp^2\right) -\left((p^2-1)(g-1)+ sp(p-1)\right)\nonumber\\
			&=\ (p^2+p+1)(g-1) + sp\nonumber\\
			&\geq\ 3.
		\end{align*}
		Thus we have $\dim(\mc{M}_L^{\boldsymbol{\alpha},rs})-\dim(Z_{\boldsymbol{\alpha}})\,\geq\, 3$,
		which completes the proof (see \eqref{es}).
	\end{proof}
	
	\section{Brauer group of moduli stack of parabolic symplectic bundles}\label{section:brauer-group-of-stack}
	
	In this section, we continue to work with the \textit{full-flag} system of multiplicities 
	$\boldsymbol{m}$, meaning $m^i_{_{x}} \,=\, 1$ for all $x\,\in\, S$ and $i$ (see Definition 
	\ref{def:parabolic-bundles}). For an even positive integer $r$, consider parabolic 
	symplectic vector bundles $(E_*,\ \varphi_*)$ of rank $r$ and having system of weights 
	$\boldsymbol{\alpha}$ such that the symplectic form takes values in
	the line bundle $L$. In case 
	$\boldsymbol{\alpha}$ does not contain $0$, then using \cite[Lemma 3.1]{BCD24} 
	we get an induced symplectic form $\varphi$ on the underlying vector bundle $E$ which takes values 
	in $L(-S)$. Consequently, the determinant of the underlying vector bundle $E$ is fixed.
	Recall that if $(E_*,\, \varphi_*)$ is a parabolic symplectic vector bundle taking values in some
	line bundle $L$ and $\eta$ is a $2$--torsion line bundle together with an isomorphism $\rho
	\,:\, \eta^{\otimes 2}\,\stackrel{\simeq}{\longrightarrow}\,\mc{O}_X$, then $(E_*\otimes \eta, \,
	\varphi_*\otimes\rho)$ is another parabolic symplectic vector bundle taking values in the same
	line bundle $L$. This leads us to the following definition.
	
	\begin{definition}\label{def:parabolic-psp-bundle}
		Fix a positive even integer $r$. A \textit{(stable) parabolic} $\psp(r,\bb{C})$--\textit{bundle} is an
		equivalence class of (stable) parabolic symplectic vector bundles taking values in a fixed line bundle $L$,
		where two parabolic symplectic vector bundles $(E_*,\ \varphi_*)$ and $(E'_*, \ \varphi'_*)$, with both
		taking values in the same line bundle $L$, are said to be \textit{equivalent} if there exists a $2$--torsion
		line bundle $\eta$ together with an isomorphism $\rho \,:\, \eta^{\otimes 2} \,
		\stackrel{\simeq}{\longrightarrow}\, \mc{O}_X$ such that there exists a parabolic isomorphism $(E'_*, \,
		\varphi'_*) \,\simeq\, (E_*\otimes \eta, \ \varphi_*\otimes \rho)$ in the sense of Definition \ref{def:parabolic-symplectic-isomorphism}.
		
	\end{definition}\begin{remark}\label{rem:psp-non-parabolic}
		If one ignores the parabolic structure, then the $\psp(r,\bb{C})$--bundles can be defined analogously
		as in Definition \ref{def:parabolic-psp-bundle}. As already remarked at the beginning of this section, if a
		parabolic $\psp(r,\bb{C})$--bundle is represented by a parabolic symplectic vector bundle $(E_*,\ \varphi_*)$ taking
		values in a line bundle $L$, then the underlying vector bundle $E$ has a symplectic form $\widetilde{\varphi}$
		which takes values in $L(-S)$, provided $0\,\notin \,\boldsymbol{\alpha}$. Consider the principal
		$\psp(r,\bb{C})$--bundle represented by the equivalence class of $(E,\ \widetilde{\varphi})$. By the underlying
		$\psp(r,\bb{C})$--bundle of a parabolic $\psp(r,\bb{C})$--bundle represented by $(E_*,\ \varphi_*)$, where
		$\varphi_*$ takes values in $L$, we shall mean the principal $\psp(r,\bb{C})$--bundle $(E,\
		\widetilde{\varphi})$, where $\widetilde{\varphi}$ takes values in $L(-S)$. 
	\end{remark}
	
	Let $J\, :=\, \begin{pmatrix}
		0 & I_p \\
		-I_p & 0
	\end{pmatrix}$, where $p\,:=\, \frac{r}{2}$, and let $\text{Gp}(r,\bb{C})$ denote the conformally
	symplectic group, meaning
	\begin{align*}
		\text{Gp}(r,\bb{C})\ :=\ \{A\in \text{GL}(r,\bb{C})\, \ \big\vert \ A^tJ A\, =\, cJ\, \ \text{for some}\ c
		\,\in\, \bb{C}^*\}.
	\end{align*}
	It can be shown that the algebraic principal $\text{Gp}(r,\bb{C})$--bundles on $X$ correspond to
	the algebraic vector bundles $E$ of rank $r$ on $X$ equipped with a non-degenerate skew-symmetric bilinear form
	$E\otimes E \,\longrightarrow\, L'$ for some line bundle $L'$ on $X$. The group $\text{Gp}(r,\bb{C})$ has
	center $\bb{C}^*$, and it fits into the short exact sequence
	\begin{align*}
		1\,\longrightarrow\, \bb{C}^*\,\longrightarrow\,\text{Gp}(r,\bb{C})\,\longrightarrow\,
		\text{PSp}(r,\bb{C})\,\longrightarrow\, 1.
	\end{align*}
	
	The associated long exact sequence of cohomologies gives the following:
	\begin{equation}\label{ee}
		H^1(X,\,\mc{O}_X^*) \ \longrightarrow\ H^1(X,\,\text{Gp}(r,\bb{C}))\ \xrightarrow{\,\,\,\delta\,\,\,}\ H^1(X,
		\, \psp(r,\bb{C})) \ \longrightarrow\ H^2(X,\,\mc{O}_X^*).
	\end{equation}
	Since $H^2(X,\, \mc{O}_X)\,=\, 0\,=\, H^3(X,\, 2\pi\sqrt{-1}{\mathbb Z})$, from the exact sequence of
	cohomologies
	$$
	H^2(X,\, \mc{O}_X)\ \longrightarrow\ H^2(X,\,\mc{O}_X^*) \ \longrightarrow\ H^3(X,\, 2\pi\sqrt{-1}{\mathbb Z})
	$$
	associated to the exponential sequence
	$$
	0 \ \longrightarrow\ 2\pi\sqrt{-1}{\mathbb Z}\ \longrightarrow\ \mc{O}_X\ \xrightarrow{\,\,\, \exp\,\,\,}\
	\mc{O}_X^* \ \longrightarrow\ 0
	$$
	it follows that $H^2(X,\,\mc{O}_X^*)\,=\,0$. Consequently, the homomorphism $\delta$ in \eqref{ee} is surjective.
	Thus the underlying principal $\psp(r,\bb{C})$--bundle of a parabolic $\psp(r,\bb{C})$--bundle (see Remark \ref{rem:psp-non-parabolic}) can be represented
	by a class of a principal $\text{Gp}(r,\bb{C})$--bundle, namely a symplectic vector bundle $(E,\,\varphi)$
	taking values in some line bundle $\mc{L}$.
	
	Let $\mf{N}^{\boldsymbol{\alpha},d}_L$
	denote the moduli stack of full-flag type stable parabolic $\psp(r,\bb{C})$--bundles $[E_*]$ of topological type $d\in\{0,1\}$, where $d\equiv \deg(L)\ \text{(mod 2)}$.
	Here, full-flag type means that its equivalence class in Definition \ref{def:parabolic-psp-bundle} can be represented by
	a stable parabolic symplectic vector bundle with full-flag systems of multiplicities (see Definition \ref{def:parabolic-bundles}); the symplectic
	form takes values in a line bundle $L$ on $X$. Let $\mc{N}^{\boldsymbol{\alpha},d}_L$ denote the coarse moduli space of $\mf{N}^{\boldsymbol{\alpha},d}_L$ . It is clear from the above description that $\mf{N}^{\boldsymbol{\alpha},d}_L$ is the quotient stack $[\mc{M}^{\boldsymbol{\alpha}}_L/\Gamma]$\ , while $\mc{N}^{\boldsymbol{\alpha},d}_L$ is the quotient variety $\mc{M}^{\boldsymbol{\alpha}}_L/\Gamma$. 
	
	\begin{remark}
		One can also define the moduli stack and the coarse moduli space of stable parabolic $\psp(r,
		\bb{C})$--bundles for partial flags as well. The case of partial flags will be considered
		in Section \ref{section:concentrated-weights}.
	\end{remark} 
	
	\begin{theorem}\label{thm:brauer-group-of-moduli-stack}
		Let $\boldsymbol{\alpha}$ be a system of weights corresponding to full-flag systems of multiplicities at each
		parabolic point (see Definition \ref{def:parabolic-bundles}). Let $\left(\mc{N}^{\boldsymbol{\alpha},
			d}_L\right)^{sm}$ denote the smooth locus of $\mc{N}^{\boldsymbol{\alpha},d}_L$. Then,
		$$\Br\left(\mf{N}^{\boldsymbol{\alpha},d}_L\right)\,\ \simeq\,\
		\Br\left(\left(\mc{N}^{\boldsymbol{\alpha},d}_L\right)^{sm}\right).$$
	\end{theorem}
	
	\begin{proof}
		Let $Z_{\boldsymbol{\alpha}}$ be as in Corollary \ref{cor:free-action-complement-codimension}, and denote
		$\mc{U}\,:=\,\mc{M}^{\boldsymbol{\alpha},rs}_L\setminus Z_{\boldsymbol{\alpha}}$\ . Consider the following diagram:
		\begin{align}\label{diagram:diag-4}
			\xymatrix{[\mc{U}/\Gamma]\ \ar@{^{(}->}[r] \ar[d]\ &\ [\mc{M}^{\boldsymbol{\alpha},rs}_L/\Gamma]\ \ar@{^{(}->}[r] \ar[d] \ &\ \mf{N}^{\boldsymbol{\alpha},d}_L \ar[d] \ \ar@{=}[r]\ & \ [\mc{M}^{\boldsymbol{\alpha}}_L/\Gamma]\\
				\mc{U}/\Gamma \ \ar@{^{(}->}[r]\ &\ \mc{M}^{\boldsymbol{\alpha},rs}_L/\Gamma \ \ar@{^{(}->}[r]\ &\ \mc{N}^{\boldsymbol{\alpha},d}_L \ \ar@{=}[r]\ & \ \mc{M}^{\boldsymbol{\alpha}}_L/\Gamma
			}
		\end{align}
		where each of the top horizontal arrows correspond to inclusions of open sub-stacks, while
		all the bottom horizontal arrows correspond to inclusions of open sub-schemes.
		Taking quotient by an action of the finite group $\Gamma$ is a finite morphism, and
		hence the quotient map preserves codimension. Therefore, it follows from Corollary \ref{cor:free-action-complement-codimension} that the complement of $\mc{U}/\Gamma$ in $\left(\mc{N}^{\boldsymbol{\alpha},d}_L\right)^{sm}$ is of codimension at least $3$. For a similar reason, the complement of the open sub-stack $[\mc{U}/\Gamma]$ in $\mf{N}^{\boldsymbol{\alpha},d}_L$ is of codimension at least $3$ as well. As $\mf{N}^{\boldsymbol{\alpha},d}_L$ is a Deligne--Mumford stack, it follows that
		\begin{equation*}
			\Br\left([\mc{U}/\Gamma]\right)\,\ \simeq\, \ \Br\left(\mf{N}^{\boldsymbol{\alpha},d}_L\right)
		\end{equation*}
		(see \cite[Proposition 4.2]{BCD23}). Now, as the action of $\Gamma$ on $\mc{U}$ is free, the left-most
		vertical arrow in the diagram \eqref{diagram:diag-4} is an isomorphism. Also, it is well--known that
		$\mc{M}^{\boldsymbol{\alpha},rs}_L$ is precisely the smooth locus of $\mc{M}^{\boldsymbol{\alpha}}_L$, and
		thus $\mc{U}$ is smooth. As $\Gamma$ acts freely on $\mc{U}$, the quotient $\mc{U}/\Gamma$ is also
		smooth. The complement of $\mc{U}/\Gamma$ in $\left(\mc{N}^{\boldsymbol{\alpha},d}_L\right)^{sm}$
		clearly has codimension at least $3$ as well. Thus we have
		$$\Br\left(\mf{N}^{\boldsymbol{\alpha},d}_L\right)\ \simeq \ Br\left([\mc{U}/\Gamma]\right)\ \simeq\
		\Br\left(\mc{U}/\Gamma\right)\ \simeq\ \Br\left(\left(\mc{N}^{\boldsymbol{\alpha},d}_L\right)^{sm}\right),$$ 
		where the last isomorphism follows from \textnormal{\cite[Theorem 1.1]{Ce19}}. This
		completes the proof.
	\end{proof}
	
	\section{Fixed-point locus of the non-parabolic symplectic moduli}\label{section:non-parabolic-case}
	
	In this section, versions of Proposition 
	\ref{prop:fixed-point-morphism-surjective} and Corollary \ref{cor:free-action-complement-codimension}
	are proved for the moduli space of usual (non-parabolic) semistable symplectic vector bundles on a curve.
	These will be used in Section \ref{section:concentrated-weights} in the computation of the Brauer
	group of the parabolic symplectic moduli. It should be clarified that even though the results in this section
	are similar to Proposition \ref{prop:fixed-point-morphism-surjective} and Corollary
	\ref{cor:free-action-complement-codimension}, their proofs crucially used the condition of \textit{full-flag}
	systems of multiplicities, and hence they can't be directly applied to the case of usual (non-parabolic) 
	symplectic vector bundles.
	
	Let $\mc{M}_L$ denote the moduli space of semistable symplectic vector
	bundles $(E,\, \varphi)$ on $X$ of rank $r$ 
	such that the symplectic form takes values in the line bundle $L$. As before, the group of $2$--torsion line bundles on $X$ act on
	$\mc{M}_L$ by tensor product. To describe this action, take a nontrivial line bundle $\eta$ on $X$ of order two,
	and fix an isomorphism
	$$
	\rho \ :\ \eta^{\otimes 2}\ \stackrel{\simeq}{\longrightarrow}\ \mc{O}_X
	$$
	as in \eqref{rho}. The action of
	$\eta$ sends a symplectic vector bundle $(E,\, \varphi)$ to $(E\otimes \eta ,\, \varphi\otimes \rho)$. Thus, a
	fixed point under the action of $\eta$ is a symplectic vector bundle $(E ,\, \varphi)$ together with an
	isomorphism of symplectic vector bundles on $X$
	$$\theta \ :\ ( E , \, \varphi)\ \stackrel{\simeq}{\longrightarrow}\ (E\otimes \eta ,\, \varphi\otimes\rho)$$
	such that the diagram
	\begin{align}\label{diagram:fixed-point-2}
		\xymatrix{E\otimes E \ar[rr]^{\varphi} \ar[d]_{\theta\otimes \theta} && L \\
			E\otimes E\otimes \eta^2\ar[urr]_{\varphi\otimes \rho} &&
		}
	\end{align}
	is commutative. Next, we note that most of the discussions in Section \ref{section:fixed-point-locus} remain valid for any 
	parabolic structure; in particular, it is applicable to the usual (non-parabolic) symplectic vector bundles 
	as well, by treating them as parabolic vector bundles with the trivial parabolic structure (see Remark 
	\ref{rem:special-structure}). Thus, starting with a vector bundle $F$ on $Y$ together with a skew-symmetric 
	bilinear form $$\phi\,:\,F\otimes F\,\longrightarrow \,\gamma^*L ,$$ we can choose a nowhere vanishing 
	section of $\eta$ and construct as before the associated spectral curve $\gamma\,:\,Y\,\longrightarrow\,
	X$, and construct a homomorphism exactly similar to \eqref{eqn:eq-1}:
	\begin{align}\label{eqn:eq-1-2}
		\phi'\ :\ (\gamma_* F)\otimes (\gamma_*F)\ \longrightarrow\ L
	\end{align}
	which is again a skew-symmetric bilinear form on $\gamma_* F$.
	
	\begin{lemma}\label{lem:non-degenerate-2}
		For any semistable symplectic vector bundle $(F,\ \phi)$ on $Y$ taking values in $\gamma^*L$, the
		direct image $(\gamma_*F,\, \phi')$ in \eqref{eqn:eq-1-2} is a semistable symplectic vector bundle
		on $X$ taking values in $L$.
	\end{lemma}
	
	\begin{proof}
		The proof for the most part is exactly same as in Lemma \ref{lem:non-degenerate} applied to the trivial 
		parabolic structure. The only modification required is in 
		showing the semistability of $(\gamma_* F,\, \phi')$, because the proof of Lemma \ref{lem:non-degenerate} uses 
		\cite[Lemma 3.3]{BCD23} which assumes full-flag systems of multiplicities.
		However, note that the condition that $(F,\, \phi)$ is
		semistable symplectic implies that the vector bundle $F$ is semistable \cite{S08}. It follows that
		$\gamma_* F$ is also semistable 
		\cite[Proposition 3.1 (ii)]{NR75}, and thus $\left(\gamma_*F,\, \phi'\right)$ is a semistable symplectic 
		vector bundle, again by \cite{S08}.
	\end{proof}
	
	\begin{lemma}\label{lem:fixed-point-2}
		Fix an isomorphism $\rho\,:\, \eta^{\otimes 2}\,\stackrel{\simeq}{\longrightarrow}\,
		\mc{O}_X$. The semistable symplectic vector bundle
		$(\gamma_*F, \, \phi')$ constructed in Lemma \ref{lem:non-degenerate-2} is a fixed point for the action of the line bundle $\eta$ on the moduli space $\mc{M}_L$.
	\end{lemma}
	
	\begin{proof}
		This is Lemma \ref{lem:fixed-point} applied to the
		trivial parabolic structure. The same proof goes through.
	\end{proof}
	
	Next, consider the other construction through which symplectic vector bundles on $X$ can arise from
	the spectral curve $Y$. As in Section \ref{section:fixed-point-locus}, let $F$ be a vector bundle on $Y$
	equipped with a pairing $\psi \, :\, F\otimes \sigma^* F \, \longrightarrow\, \gamma^*L$
	which is non-degenerate in the sense that the homomorphism
	\begin{equation}\label{ei-2}
		F\, \longrightarrow\, \gamma^*L\otimes (\sigma^* F)^\vee\,=\, \gamma^*L\otimes \sigma^* (F^\vee)
	\end{equation}
	induced by $\psi$ is an isomorphism. Also, consider
	$$
	\sigma^*\psi\, :\, \sigma^*(F\otimes \sigma^* F)\,=\, (\sigma^* F)\otimes F \, \longrightarrow\,
	\sigma^*\gamma^*L\,=\, \gamma^*L.
	$$
	As before, assume that 
	\begin{align}\label{eqn:extra-condition-2}
		\sigma^*\psi\,=\, -\psi\circ \varpi,
	\end{align}
	where $\varpi\, :\, (\sigma^* F)\otimes F\, \longrightarrow\, F\otimes \sigma^* F$ is the natural involution that switches
	the two factors of the tensor product. 
	Clearly, we have an analogue of \eqref{e5} in the usual (non-parabolic) case, which is denoted by
	\begin{equation*}
		\mathcal{J}_{_{F}}\ : \ \gamma_* (F)\otimes \gamma_* (F)\ \longrightarrow \gamma_* (F \otimes \sigma^* F).
	\end{equation*} 
	Next, consider the following composition of homomorphisms:
	$$
	(\gamma_* F)\otimes (\gamma_* F)\ \xrightarrow{\,\,\, \mathcal{J}_{_{F}}\,\,\,} \ \gamma_* (F \otimes \sigma^* F)
	\ \xrightarrow{\gamma_*(\psi)}\ \gamma_*\gamma^* L\ \longrightarrow\ L\ .
	$$
	Denote this composition of homomorphisms by
	\begin{equation}\label{e6-2}
		\psi'\ :\ (\gamma_* F)\otimes (\gamma_* F)\ \longrightarrow \ L.
	\end{equation}
	
	\begin{lemma}\label{lem:another-fixed-point-2}
		Let $F$ be a semistable vector bundle on $Y$ together with a non-degenerate pairing
		$$\psi\ :\ F\otimes (\sigma^*F)\ \longrightarrow\ \gamma^*L$$ satisfying \eqref{eqn:extra-condition-2}.
		Then the following two statements hold:
		\begin{enumerate}[(1)]
			\item $(\gamma_*F, \ \psi')$ in in \eqref{e6-2}
			is a semistable symplectic vector bundle on $X$ taking values in the line bundle $L$.
			
			\item Fix an isomorphism $\rho\,:\, \eta^{\otimes 2}\,\stackrel{\simeq}{\longrightarrow}\, \mc{O}_X$. The
			semistable symplectic vector bundle
			$(\gamma_*(F), \, \psi')$ in (1) is a fixed point for the action of the line bundle $\eta$ on the moduli space $\mc{M}_L$.
		\end{enumerate}
	\end{lemma}
	
	\begin{proof}
		This is Lemma \ref{lem:another-fixed-point} applied to the trivial parabolic structure.
	\end{proof}
	
	\vspace{0.2cm}
	
	\subsection{Codimension estimation of the fixed point locus for non-parabolic case}\hfill\\
	
	Fix an even positive integer $r$ and also fix a line bundle $L$ on $X$. Recall that $\mc{M}_L$ denotes
	the moduli space of semistable symplectic vector bundles $(E,\, \varphi)$ on $X$ of rank $r$
	such that the symplectic form takes values in $L$. Let 
	$$\mc{M}^{rs}_L\ \,\subset\ \,\mc{M}_L$$ denote the 
	open Zariski subset consisting of regularly stable symplectic vector bundles. It is 
	straightforward to see that $\mc{M}^{rs}_L$ is invariant under the 
	action on $\mc{M}_L$ of the $2$--torsion line bundles on $X$.
	
	Next, the counterpart of Proposition \ref{prop:fixed-point-morphism-surjective} will be proved in
	the non-parabolic set-up. Some changes are necessary in the non-parabolic setting, which
	will be described below. Let $\mc{N}_{Y}$ denote the moduli space of semistable symplectic vector bundles
	of rank $\frac{r}{2}$ on $Y$ such that the symplectic form takes values in $\gamma^*L$. Note that
	$\mc{N}_{Y}$ is empty if $\frac{r}{2}$ is odd. Also, let $\mc{N}'_Y$ denote the moduli space of semistable
	vector bundles $F$ of rank $\frac{r}{2}$ on $Y$ equipped with a pairing $\psi\, :\, 
	F\otimes\sigma^*F \,\longrightarrow \,\gamma^*L$ that satisfies \eqref{eqn:extra-condition-2}.
	
	Let $\left(\mc{M}_L\right)^{\eta}\, \subset\, \mc{M}_L$ denote the fixed-point locus
	for the action of $\eta$ on $\mc{M}_L$.
	By Lemma \ref{lem:fixed-point-2} and Lemma \ref{lem:another-fixed-point-2}, there exists a
	natural morphism \begin{align}\label{eqn:morphism-to-fixed-point-2}
		f\ :\ \mc{N}_{Y}\coprod \mc{N}'_Y\ \longrightarrow\ 
		\left(\mc{M}_L\right)^{\eta}.
	\end{align}
	
	\begin{proposition}\label{prop:fixed-point-morphism-surjective-2}
		Let $V \,:=\, f^{-1}\left((\mc{M}_L^{rs})^{\eta}\right)$ with $f$
		being the map in \eqref{eqn:morphism-to-fixed-point-2}. The restricted morphism
		$$f\big\vert_V\ :\ V\ \longrightarrow\ (\mc{M}_{L}^{rs})^{\eta}$$ is surjective.
	\end{proposition}
	
	\begin{proof}
		The proof is essentially the same as that of Proposition \ref{prop:fixed-point-morphism-surjective}
		applied to the trivial parabolic structure. The only change required here is in the fact that in Proposition \ref{prop:fixed-point-morphism-surjective}, we used
		the proof in \cite[Lemma 3.3]{BCD23} to produce a parabolic symplectic vector bundle on $Y$; while here,
		we need to use the proof in \cite[Lemma 2.1]{BH10} to conclude the same.
		The rest of the proof remains exactly the same.
	\end{proof}
	
	\begin{corollary}\label{cor:free-action-complement-codimension-2}
		Let $\Gamma$ denote the group of $2$--torsion line bundles on $X$. Consider the
		Zariski closed subset 
		$$Z\ :=\ \bigcup_{\eta\in\Gamma\setminus\{\mc{O}_X\}}
		(\mc{M}_L^{rs})^{\eta}\ \subset\ \mc{M}^{rs}_L .$$
		The codimension of $Z$ in $\mc{M}^{rs}_L$
		is at least $3$.
	\end{corollary}
	
	\begin{proof}
		The argument is almost similar to that of Corollary \ref{cor:free-action-complement-codimension}. Let
		$\eta\,\in\, \Gamma$ be a nontrivial line bundle. Let $\gamma\,:\,Y\,\longrightarrow\,
		X$ be the spectral curve corresponding to $\eta$. To prove the result it is enough
		to show that
		\begin{equation}\label{es-2}
			\dim(\mc{M}^{rs}_L)-
			\dim((\mc{M}^{rs}_L)^{\eta})\,\geq\, 3.
		\end{equation}
		Now, by Proposition
		\ref{prop:fixed-point-morphism-surjective-2},
		$$\dim\left((\mc{M}_L^{rs})^{\eta}\right)\ \leq\
		\max\left\{\dim(\mc{N}_Y),\ \dim(\mc{N}'_Y)\right\}.$$
		Let $g(Y)$ denote the genus of $Y$. By Riemann--Hurwitz formula, $g(Y) \,=
		\, 2(g-1)+1$. For notational convenience, denote $p\,:=\, \frac{r}{2}$. We have the
		following expressions for dimensions:
		\begin{align*}
			&\dim(\mc{M}^{rs}_L) \,=\, \dim(\mc{M}_L)\,=\,
			p(2p+1)(g-1)\ \ \ \ \textnormal{\cite[Lemma 3.10]{BR88}},\\
			\text{while}\ \ &\dim(\mc{N}_Y) \,=\, 0 \ \ \textnormal{(if}\ \ p \ \ \textnormal{is odd), and} \nonumber\\
			& \dim(\mc{N}_Y)\, =\,\dfrac{p}{2}(p+1) \left(2g(Y)-1\right)
			\ \ \textnormal{(if}\ p \ \textnormal{is even)}\nonumber\\
			& \qquad \qquad =\ \frac{p}{2}\left(p+1\right)(2g-2)\nonumber\\
			&\qquad\qquad =\ p(p+1)(g-1).
		\end{align*}
		For a vector bundle $F$ in $\mc{N}'_Y$, the line bundle
		$\det(\gamma_*F)$ is fixed. Thus, we get the following:
		\begin{align*}
			\dim(\mc{N}'_Y) &\ \leq\ p^2\left(g(Y) -1\right) +1 -g \nonumber\\
			&=\ p^2(2g-2) +1 -g \nonumber\\
			&=\ (p^2-1)(g-1).
		\end{align*}
		It follows that 
		\begin{align*}
			\dim(\mc{M}^{rs}_L) - \dim(\mc{N}_Y)\, &
			=\, \left(p(2p+1)(g-1)\right) - \left(p(p+1)(g-1)\right)\nonumber\\
			& =\,p^2(g-1)\nonumber\\
			&\geq\ 3 .
		\end{align*}
		\begin{align*}
			\text{Similarly},\, \ \dim(\mc{M}^{rs}_L)-
			\dim(\mc{N}'_Y) \,&\geq\, \left(p(2p+1)(g-1)\right) -\left((p^2-1)(g-1)\right)\nonumber\\
			&=\ (p^2+p+1)(g-1)\nonumber\\
			&\geq\ 3.
		\end{align*}
		Thus we have $\dim(\mc{M}_L^{rs})-\dim Z\,\geq\, 3$,
		which completes the proof (see \eqref{es-2}).
	\end{proof}
	
	\section{Brauer groups for concentrated weights}\label{section:concentrated-weights}
	
	Henceforth, the system of multiplicities are not needed to be of full-flag type.
	
	As it was observed in \cite{BCD24}, for the parabolic symplectic set-up 
	it is necessary to assume that the system of multiplicities $\boldsymbol{m}$ and weights 
	$\boldsymbol{\alpha}$ are of \textit{symmetric type} (see Definition \ref{def:concentrated-weights} below). We
	shall begin with a particular version of such types of weights, namely a \textit{concentrated} system
	of weights. For convenience, these two notions are recalled first.
	
	\begin{definition}[{\cite[Definition 3.4 and Definition 3.7]{BCD24}}]\label{def:concentrated-weights}
		Let $r$ be a positive even integer. Fix parabolic points $S\, \subset\, X$, and also fix a subset of positive
		integers $\{\ell(p)\}_{p\in S}$\, satisfying the condition $\ell(p)\ \leq\ r$ for all\ $p\,\in\, S$.
		Suppose that
		\begin{align*}
			\boldsymbol{m} \,\,=\,\, \left\{\left(m^1_{_{p}},\, m^2_{_{p}},\, \cdots,\, m^{\ell(p)}_{_{p}}\right)_{p\in S}\right\},\ \
			\, \boldsymbol{\alpha}\,\,=\,\,\left\{\left(\alpha^1_{_{p}}\,<\,\alpha^2_{_{p}}\,<\,
			\cdots\,<\,\alpha^{\ell(p)}_{_{p}}\right)_{p\in S}\right\}
		\end{align*}
		are respectively the systems of multiplicities and weights on points of $S$ (so $\sum_{i=1}^{\ell(p)}m^i_{_{p}}
		\,=\, r$ for all $p\,\in\, S)$.
		\begin{itemize}
			\item $\boldsymbol{m}$ is said to be of \textit{symmetric type} if $m^j_{_{p}} \,=\,
			m^{^{\ell(p)+1-j}}_{_{p}}$ for all $p\,\in\, S$ and $1\,\leq\, j\,\leq\, \ell(p)$.
			
			\item $\boldsymbol{\alpha}$ is said to be of \textit{symmetric type} if $\alpha^j_{_{p}}
			\,=\, 1-\alpha^{^{\ell(p)+1-j}}_{_{p}}$ for all $p\,\in\, S$ and $1 \, \leq\, j \, \leq \, \ell(p)$.
			
			\item The system of weights $\boldsymbol{\alpha}$ is called
			\textit{concentrated} if it of symmetric type
			and satisfies the inequality $\underset{p\in S}{\sum} \left(\frac{1}{2}-\alpha^1_{_{p}}\right) \, <\,
			\dfrac{1}{r^2}$.
		\end{itemize}
	\end{definition}
	
	Fix an even positive integer $r$, a set of parabolic points $S$ on $X$ and a system of multiplicities
	$\boldsymbol{m}$ of symmetric type. Let $\boldsymbol{\alpha}$ be a concentrated system of weights
	(Definition \ref{def:concentrated-weights}). Denote $L(S)\, :=\, L\otimes \mc{O}_X(S)$. Recall from
	Section \ref{section:fixed-point-locus} that the group $\Gamma$ of $2$--torsion line bundles on $X$ acts
	on the semistable moduli space $\mc{M}^{\boldsymbol{m,\alpha}}_{L(S)}$ through tensorization. For the sake of
	convenience, the following definition is recalled.
	
	\begin{definition}\label{def:psp-moduli}
		Fix a line bundle $\mc{L}$ on the curve $X$, and take $d\,\in\,\{0,\,1\}$
		such that $\deg(\mc{L}) \,\equiv\, d \ (\text{mod}\ 2)$. The (twisted) coarse moduli
		space of semistable parabolic $\psp(r,\bb{C})$--bundles on $X$ of topological type $d$, which can be represented by stable
		parabolic symplectic vector bundles $(E_*,\, \varphi_*)$ with $\varphi_*$ taking values in $\mc{L}$ (see 
		Definition \ref{def:parabolic-psp-bundle}), is
		the quotient variety
		$$\mc{N}^{\boldsymbol{m,\alpha},d}_{\mc{L}}\ :=\ \mc{M}^{\boldsymbol{m,\alpha}}_{\mc{L}}/\Gamma .$$ 
	\end{definition}
	
	Fix a concentrated system of weights $\boldsymbol{\alpha}$; we assume that $\boldsymbol{\alpha}$ does not contain $0$.
	As before, let $\mc{M}_L$ denote the coarse moduli space of semistable symplectic vector bundles of
	rank $r$ on $X$ with the symplectic form taking values in a line bundle $L$. Let $\mc{M}^{rs}_L\,\subset\,
	\mc{M}_L$ denote the Zariski open subvariety of regularly stable parabolic symplectic vector
	bundles. It follows from \cite[Lemma 4.1]{BCD24} that there exists a morphism
	\begin{align*}
		\pi_0\ :\ \mc{M}^{\boldsymbol{m,\alpha}}_{L(S)}\ \longrightarrow\ \mc{M}_L,
	\end{align*} 
	whose restriction to $\pi_0^{-1}\left(\mc{M}^{rs}_L\right)$ (which we also denote by $\pi_0$ by a slight abuse of notation), namely
	\begin{align}\label{eqn:fibr}
		\pi_0\ :\ \pi_0^{-1}\left(\mc{M}^{rs}_L\right)\ \longrightarrow\ \mc{M}^{rs}_L,
	\end{align}
	is a fiber bundle map, with fibers isomorphic to the isotropic flag variety
	\begin{equation}\label{ef}
		F \ :=\ \prod_{i=1}^{|S|} \text{Sp}(r,\bb{C})/P_i,
	\end{equation}
	where $P_i\, \subset\, \text{Sp}(r,\bb{C})$ is the parabolic
	subgroup corresponding to the flag at the $i$--th parabolic point \cite[Lemma 4.1]{BCD24}. The morphism $\pi_0$ in \eqref{eqn:fibr} is clearly equivariant for the actions of $\Gamma$ on
	$\mc{M}^{\boldsymbol{m,\alpha}}_{L(S)}$ and $\mc{M}_L$. Evidently, $\mc{M}^{rs}_L$ is a $\Gamma$--invariant subvariety. Let 
	\begin{align}\label{eqn:V}
		V\ :=\ \mc{M}^{rs}_L\setminus Z\ \subset\ \mc{M}^{rs}_L,
	\end{align}
	where $Z$ is defined in Corollary
	\ref{cor:free-action-complement-codimension-2}. The open subset $V$ in \eqref{eqn:V}
	is easily seen to be $\Gamma$--invariant. It follows that $U\,=\,\pi_0^{-1}(V)$ is also $\Gamma$--invariant, and the map $\pi$ in \eqref{eqn:fibration} is $\Gamma$--equivariant. 
	
	\begin{lemma}\label{lem:picard-Brauer-V}
		Consider the $\Gamma$--invariant subvariety $V$ as in \eqref{eqn:V}.
		The following holds: $$\Br(V)\ \stackrel{\simeq}{\longrightarrow} \ \Br(\mc{M}^{rs}_L).$$
		
		The group $\Gamma$ acts trivially on $\Pic(V)$.
	\end{lemma}
	
	\begin{proof}
		The isomorphism $\Br(V)\ \stackrel{\simeq}{\longrightarrow} \ \Br(\mc{M}^{rs}_L)$ follows immediately
		from \cite[Theorem 1.1]{Ce19}, the codimension estimate in Corollary 
		\ref{cor:free-action-complement-codimension-2} and the fact that $\mc{M}^{rs}_L$ is smooth
		\cite[Proposition 2.3]{BH12}.
		
		To see that $\Gamma$ acts trivially on $\Pic(V)$, first consider the cases where $\mc{M}_L$ is locally factorial, 
		which holds in the following two cases (see \cite[Theorem (1.6), p. 501]{LS97} and \cite[Corollary 8.2]{BHol13}):
		\begin{enumerate}
			\item $d\,=\,0$ (equivalently, $\deg(L)$ is even),
			
			\item $d\,=\,1$ (equivalently, $\deg(L)$ is odd) and $\frac{r}{2}$ is odd.
		\end{enumerate}
		As the variety $\mc{M}_L$ is normal and $\mc{M}^{rs}_L$, is precisely the smooth locus of $\mc{M}_L$
		\cite[Proposition 2.3]{BH12}, it follows that the complement of $\mc{M}^{rs}_L$ in $\mc{M}_L$ is of codimension at least $2$. This fact, combined with the codimension estimate in Corollary \ref{cor:free-action-complement-codimension-2}, implies that the complement of $V$ in $\mc{M}_L$ is of codimension at least $2$ as well. As $\mc{M}_L$ is locally factorial in the two cases mentioned earlier, we conclude that
		\begin{align*}
			\Pic(V)\ \stackrel{\simeq}{\longrightarrow}\ \Pic(\mc{M}_L).
		\end{align*}
		Now, it is known that $\Pic(\mc{M}_L)$ is infinite cyclic \cite[1.6]{LS97}. As
		the action of $\Gamma$ must fix the ample generator of $\Pic(\mc{M}_L)$, it now
		follows that $\Gamma$ acts trivially on $\Pic(\mc{M}_L)\,\simeq\, \Pic(V)$.
		
		In the remaining case, meaning $d=1$ (equivalently, $\deg(L)$ is odd) and $\frac{r}{2}$ is even, we have an 
		inclusion $\Pic(\mc{M}_L)\ \hookrightarrow\ \Pic(V)$. The Picard group $\Pic(\mc{M}_L)$
		(respectively, $\Pic(V)$) is infinite cyclic, and it is generated by the smallest power of the generator
		of the Picard group of the affine Grassmannian that descends to $\mc{M}_L$ (respectively, $V$) (see 
		\cite{BLS98}). It follows that the inclusion $\Pic(\mc{M}_L)\,\hookrightarrow\,\Pic(V)$ is of the form $\mc{L}\,
		\longmapsto\, \mc{L}^d$ for some positive integer $d$, where $\mc{L}$ is the generator of $\Pic(\mc{M}_L)$. It has been 
		already argued that $\Gamma$ acts trivially on $\Pic(\mc{M}_L)$. It now follows immediately that $\Gamma$ acts 
		trivially on $\Pic(V)$ as well.
	\end{proof}
	
	Let $V$ be as in \eqref{eqn:V}. Consider
	\begin{equation}\label{du}
		U\ :=\ \pi_0^{-1}(V)\ \subset\ \pi_0^{-1}\left(\mc{M}^{rs}_L\right),
	\end{equation}
	and denote the restriction of $\pi_0$ to $U$ by $\pi$, namely
	\begin{align}\label{eqn:fibration}
		\pi \ :\ U\ \longrightarrow\ V .
	\end{align}
	This yields the following commutative diagram:
	\begin{align*}
		\xymatrix{
			{U\,} \ar@{^{(}->}[r] \ar[d]_{\pi} & \pi_0^{-1}\left(\mc{M}^{rs}_{L}\right) \ar[d]^{\pi_0}\\
			{V\,} \ar@{^{(}->}[r] & \mc{M}^{rs}_L
		}
	\end{align*}
	
	\begin{lemma}\label{lem:codimension-estimate-U}
		The open subset $U$ in \eqref{du} is smooth, and its complement in the smooth locus
		$\left(\mc{M}^{\boldsymbol{m,\alpha}}_{L(S)}
		\right)^{sm}$ has codimension at least 2.
	\end{lemma}
	
	\begin{proof}
		
		It is known that $\mc{M}^{rs}_L$ is the smooth locus of $\mc{M}_L$ \cite[Proposition 2.3]{BH12}. Thus $V$ is 
		smooth. As the fibers $F$ of $\pi$ (see \eqref{ef}) are smooth rational projective varieties, and $\pi$ is a 
		fiber bundle map, it follows that $U$ is also smooth.
		
		To prove that the codimension of the complement of $U$ in $\left(\mc{M}^{\boldsymbol{m,\alpha}}_{L(S)}
		\right)^{sm}$ is at least 2, first note that since the map $\pi_0$ in \eqref{eqn:fibr} is a
		fiber bundle map, it follows from Corollary \ref{cor:free-action-complement-codimension-2} that
		$\pi_0^{-1}(Z)$ is a Zariski closed subset of codimension at least 2 in $\pi_0^{-1}\left(\mc{M}_L^{rs}\right)$.
		Evidently, $\pi_0^{-1}(Z)$ is precisely the complement of $U$ in $\pi_0^{-1}\left(\mc{M}^{rs}_L\right)$.
		Thus, from the following chain of inclusions of open subsets
		\begin{align*}
			U \ \subset \ \pi_0^{-1}\left(\mc{M}^{rs}_L\right)\ \subset\ 
			\left(\mc{M}^{\boldsymbol{m,\alpha}}_{L(S)}\right)^{sm}
		\end{align*}
		it follows easily that the complement of $U$ in $\left(\mc{M}^{\boldsymbol{m,\alpha}}_{L(S)}\right)^{sm}$
		is of codimension at least 2.
	\end{proof}
	
	\begin{lemma}\label{lem:Gamma-action-on-Pic(U)}
		The Picard group $\Pic(U)$ is torsion-free. The action of $\Gamma$ on $\Pic(U)$ is the trivial one.
	\end{lemma}
	
	\begin{proof}
		Recall the fiber bundle $\pi$ in \eqref{eqn:fibration} with fiber $F$ (see \eqref{ef}). Since $F$ is a
		smooth projective variety, and $\pi$ is a fiber bundle map, it follows immediately that $U$ is also
		smooth, as well as $H^0(F,\,\mc{O}_F^*)\,=\, \bb{C}^*$. Using \cite[Proposition 2.3]{FI73}, one obtains the
		following exact sequence of Picard groups:
		\begin{align}\label{eqn:picard-group-exact-sequence}
			0\ \longrightarrow\ \Pic(V)\ \stackrel{\pi^*}{\longrightarrow}\ \Pic(U)\ \stackrel{\omega}{\longrightarrow}
			\ \Pic(F),
		\end{align}
		where the homomorphism $\omega$ sends a line bundle on $U$ to its restriction to a fiber of $\pi$.
		
		Note that $\Pic(F)$ is a free abelian group of finite rank. In particular, it is torsionfree. Also,
		$\Pic(V)$ is torsionfree, in fact, it is isomorphic to $\mathbb Z$. Therefore, from
		\eqref{eqn:picard-group-exact-sequence} it follows that $\Pic(U)$ is torsionfree.
		
		To prove that the action of $\Gamma$ on $\Pic(U)$ is the trivial one, first note
		that the homomorphism $\pi^*$ in \eqref{eqn:picard-group-exact-sequence} is $\Gamma$--equivariant,
		because the map $\pi$ is $\Gamma$--equivariant. Therefore, the action of $\Gamma$ on $\Pic(U)$ induces
		an action of $\Gamma$ on the image $\omega(\Pic(U))$, where $\omega$ is the homomorphism in
		\eqref{eqn:picard-group-exact-sequence}.
		
		As noted before, $\Pic(F)$ is a free abelian group of finite rank. From the
		description of $\Pic(U)$ (see \cite{LS97})
		it follows that there is a subgroup ${\mathbb S}\, \subset\, \Pic(U)$ such that the following statements hold:
		\begin{itemize}
			\item The restriction of $\omega$ (see \eqref{eqn:picard-group-exact-sequence}) to $\mathbb S$ is injective.
			
			\item The image $\omega({\mathbb S})$ is a finite index subgroup of $\Pic(F)$.
			
			\item For the action of $\Gamma$ on $\Pic(U)$, every element of $\mathbb S$ is fixed by $\Gamma$.
		\end{itemize}
		To construct $\mathbb S$, consider $F$ in \eqref{ef}. Take any $1\, \leq\, i\, \leq\, |S|$, and fix
		a $\text{PSp}(r,\bb{C})$--equivariant line bundle ${\mathcal L}_i$ on $\text{Sp}(r,\bb{C})/P_i$
		(see \eqref{ef}); note that $\text{PSp}(r,\bb{C})$ acts on $\text{Sp}(r,\bb{C})/P_i$ as left-translations
		because the center of $\text{Sp}(r,\bb{C})$ lies in $P_i$. Now ${\mathcal L}_i$ produces a line bundle on $U$
		using the quasi-parabolic flag, at the $i$-th point of $S$, of the parabolic symplectic bundles. The subgroup
		${\mathbb S}\, \subset\, \Pic(U)$ consists of these line bundles.
		As $\omega({\mathbb S})$ is a finite index subgroup of $\Pic(F)$, it follows that
		$\omega({\mathbb S})$ is a finite index subgroup of $\omega(\Pic(U))$. Since the action of $\Gamma$
		on $\mathbb S$ is the trivial one, it follows that the action of $\Gamma$ on $\omega(\Pic(U))$
		fixes $\omega({\mathbb S})$ pointwise. Consequently, $\Gamma$ acts trivially on $\omega(\Pic(U))$.
		
		In Lemma \ref{lem:picard-Brauer-V} it was shown that $\Gamma$ acts trivially on $\Pic(V)$. Since the
		action of $\Gamma$ on $\omega(\Pic(U))$ is also the trivial one, we now conclude that the action of
		$\Gamma$ on $\Pic(U)$ is the trivial one.
	\end{proof}
	
	The quotient
	$$\mc{N}^{d}_L\ =\ \mc{M}_L/\Gamma$$
	is the moduli space of semistable $\psp(2r,\bb{C})$--bundles of topological type $d\,\in\,
	\{0,\,1\}$ (see Definition \ref{def:psp-moduli} for the parabolic case). As $V$
	is smooth and the $\Gamma$-action on $V$ is free,
	it follows that $V/\Gamma$ is also smooth. Thus $V/\Gamma$ is an open subset in he smooth locus $\left(\mc{N}^{d}_L\right)^{sm}$ of $\mc{N}^d_L$. As the codimension of $\mc{M}^{rs}_L\setminus V\, \subset\, \mc{M}^{rs}_L$ is at least two
	(Corollary \ref{cor:free-action-complement-codimension-2}), a straightforward codimension estimate shows that
	\begin{align}\label{eqn:codimension-estimate-V-mod-Gamma}
		\Br(V/\Gamma)\ \stackrel{\simeq}{\longrightarrow}\ \Br\left(\left(\mc{N}^{d}_L\right)^{sm}\right).
	\end{align}
	
	As the morphism $\pi$ in \eqref{eqn:fibration} is $\Gamma$--equivariant with respect to the actions of $\Gamma$
	on $U$ and $V$, it descends to a map
	$$\overline{\pi}\ :\ U/\Gamma\ \longrightarrow\ V/\Gamma .$$ 
	Consequently, we have a commutative diagram
	\begin{align}\label{diagram:quotient}
		\xymatrix{U \ar[r]^{q} \ar[d]_{\pi} & U/\Gamma \ar[d]^{\overline{\pi}} \\
			V \ar[r]^{q'} & V/\Gamma }
	\end{align}
	where $\pi$ and $\pi'$ are fiber bundle maps with fiber $F$, while $q$ and $q'$ are the quotient maps.
	As $\Gamma$ acts freely on $V$, and $\pi$ is $\Gamma$--equivariant, we conclude that
	$\Gamma$ acts freely on $U$.
	Consequently, the quotient maps $U\stackrel{q}{\longrightarrow} U/\Gamma$ and $V\stackrel{q'}{\longrightarrow}V/\Gamma$ are finite \'etale covers.
	
	It follows from Lemma \ref{lem:codimension-estimate-U} that both $U$ and $U/
	\Gamma$ are smooth. The codimension estimate in Lemma \ref{lem:codimension-estimate-U} also tells us that
	\begin{align}\label{eqn:brauer-group-U}
		\Br\left(U\right)\,\ \stackrel{\simeq}{\longrightarrow}\,\
		\Br\left(\left(\mc{M}^{\boldsymbol{m,\alpha}}_{L(S)}\right)^{sm}\right)\ \ \text{ and }\ \ 
		\Br\left(U/\Gamma\right)\,\ \stackrel{\simeq}{\longrightarrow}\ \,
		\Br\left(\left(\mc{N}^{\boldsymbol{m,\alpha},d}_{L(S)}\right)^{sm}\right)
	\end{align}
	(see \cite[Theorem 1.1]{Ce19}).
	
	\begin{proposition}\label{prop:brauer-group-U}
		Assume that one of the following three conditions is satisfied:
		\begin{enumerate}[(a)]
			\item $d\,=\,0$ (equivalently $\deg(L)$ is even) and $m^i_{_{p}}\,=\, 1$ for some $p\,\in\, S$ and $i$;\label{c1}
			
			\item $d\,=\,1$ (equivalently $\deg(L)$ is odd) and $\frac{r}{2}\,\geq\, 3$ is odd;\label{c2}
			
			\item $d\,=\,1$ (equivalently $\deg(L)$ is odd), $\frac{r}{2}\,\geq\, 3$ is even and $m^i_{_{p}}\,=\, 1$ for
			some $p\,\in\, S$ and $i$.\label{c3}
		\end{enumerate}
		The group $\Br\left(U/\Gamma\right)$ is identified with the kernel of the homomorphism $\Br(V/\Gamma)\
		\longrightarrow \ \Br(V)$ induced from the quotient map $q'\,:\,V\,\longrightarrow \,V/\Gamma$
		in \eqref{diagram:quotient}.
	\end{proposition}
	
	\begin{proof}
		We shall closely follow the argument in \cite[Proposition 5.1]{BCD23}. Recall the isomorphism $$\Br\left(U\right)\,\ \stackrel{\simeq}{\longrightarrow}\,\
		\Br\left(\left(\mc{M}^{\boldsymbol{m,\alpha}}_{L(S)}\right)^{sm}\right)$$ from \eqref{eqn:brauer-group-U}. It follows from the description of $\Br\left(\left(\mc{M}^{\boldsymbol{m,\alpha}}_{L(S)}\right)^{sm}\right)$ in \cite[Theorem 4.5]{BCD24} that under any of the conditions \eqref{c1}, \eqref{c2} or \eqref{c3}, we have 
		$$\Br(U)=0.$$
		Moreover, it follows from Lemma \ref{lem:picard-Brauer-V} and Lemma \ref{lem:Gamma-action-on-Pic(U)} that both $\Pic(U)$ and $\Pic(V)$ are torsion--free, and the actions of $\Gamma$
		on $\Pic(U)$ and $\Pic(V)$ are trivial. This leads to the following equalities:
		\begin{align}
			\Pic(V)^{\Gamma}\ =\ \Pic(V), \ \ \ H^1(\Gamma,\,\Pic(V))
			\ =\ \Hom(\Gamma,\,\Pic(V))\ =\ 0,\label{pic1}\\
			\Pic(U)^{\Gamma}\ =\ \Pic(U), \ \ \ H^1(\Gamma,\,\Pic(U))
			\ =\ \Hom(\Gamma,\,\Pic(U))\ =\ 0.\label{pic2}
		\end{align}
		In light of these equalities, the Hochschild--Serre spectral sequences associated to the finite \'etale covers
		$U\ \stackrel{q}{\longrightarrow}\ U/\Gamma$ and $V\ \stackrel{q'}{\longrightarrow}\ V/\Gamma$
		yield the following two exact sequences (see \cite[III Theorem 2.20]{Mi}):
		\begin{align}
			0 &\, \longrightarrow\, \chi(\Gamma)\, \stackrel{f}{\longrightarrow}\, \Pic(U/\Gamma)\,
			\stackrel{q^*}{\longrightarrow}\,\Pic(U)\,\stackrel{g}{\longrightarrow}\, H^2(\Gamma,\,\bb{C}^*)
			\, \longrightarrow\, \Br(U/\Gamma)\,\longrightarrow\,\Br(U)^{\Gamma}\,=\,0, \label{eq1}\\
			0 &\, \longrightarrow\, \chi(\Gamma)\, \stackrel{f'}{\longrightarrow}\, \Pic(V/\Gamma)\,
			\stackrel{q'^*}{\longrightarrow}\,\Pic(V)\,\stackrel{g'}{\longrightarrow}\, H^2(\Gamma,\,\bb{C}^*)
			\,\longrightarrow\, \Br(V/\Gamma)\,\longrightarrow\,\Br(V)\, \longrightarrow\, 0, \label{eq2}
		\end{align} 
		where $\chi(\Gamma) \,:=\, \text{Hom}(\Gamma,\bb{C}^*)$ is the character group of $\Gamma$. Note that the map $\Br(V/\Gamma)\longrightarrow\Br(V)$ in \eqref{eq2} is surjective due to \cite[Theorem 6.3]{BHol13} combined with Lemma \ref{lem:picard-Brauer-V}.
		
		Let us denote $H\,:=\, \text{coker}(q^*)$ and	$H'\,:=\, \text{coker}(q'^*)$. We claim that $$H\ \stackrel{\simeq}{\longrightarrow}\ H'.$$ To see this, consider the diagram \eqref{diagram:quotient}, which leads to the following diagram (see \eqref{eqn:picard-group-exact-sequence}):
		\begin{align}\label{diagram:picard-diagram}
			\xymatrix{ 0 \ar[r] & \Pic(V/\Gamma) \ar[r]^{\overline{\pi}^*} \ar[d]_{q'^*} &\Pic(U/\Gamma) \ar@{->>}[r] \ar[d]_{q^*}& \coker(\overline{\pi}^*)\ar@{^{(}->}[r] \ar[d]& \Pic(F) \ar[d]_{\simeq} \\
				0 \ar[r] & \Pic(V) \ar[r]^{\pi^*} & \Pic(U) \ar@{->>}[r] & \coker(\pi^*) \ar@{^{(}->}[r]& \Pic(F) }
		\end{align}
		where the rightmost vertical map is an isomorphism by \cite[Claim 5.1]{BCD23}. It follows that the induced map $\coker(\overline{\pi}^*) \longrightarrow\coker(\pi^*)$ is an isomorphism. Using snake lemma, one immediately obtains from the diagram \eqref{diagram:picard-diagram} that $\coker(q^*)\ \stackrel{\simeq}{\longrightarrow}\coker(q'^*)$, which proves our claim.
		
		From \eqref{eq1} and \eqref{eq2} we get the following two exact sequences:
		\begin{align}
			0&\ \longrightarrow\ H\ \longrightarrow\ H^2(\Gamma,\,\bb{C}^*)
			\ \longrightarrow\ \Br(U/\Gamma)\ \longrightarrow\ 0,\label{eq3}\\
			0&\ \longrightarrow \ H' \ \longrightarrow \ H^2(\Gamma,\, \bb{C}^*)
			\ \longrightarrow\ \Br(V/\Gamma)\ \longrightarrow\ \Br(V)\ \longrightarrow\
			0\label{eq4}.
		\end{align}
		Now, as the complement of the open
		subset $V\,\subset \,\mc{M}_L$ is of codimension at least 2 and $\mc{M}_L$ is a normal projective variety, we have $H^0(V,\ \bb{G}_m)=\bb{C}^*$. For the exact same reason, the open subset $U\,\subset\, \mc{M}^{\boldsymbol{m,\alpha}}_{L(S)}$ also satisfies $H^0(U,\ \bb{G}_m)=\bb{C}^*$. It follows that $\pi\,:\,U\,\longrightarrow\, V$ induces an isomorphism
		$$\bb{C}^*\ =\ H^0(V,\,\mathbb{G}_m)\ \longrightarrow\
		H^0(U,\,\mathbb{G}_m)\ =\ \bb{C}^* .$$
		This in turn produces an isomorphism $H^2(\Gamma,\,\bb{C}^*)\,\longrightarrow \,H^2(\Gamma,\,\bb{C}^*)$, which
		takes $H$ to $H'$.
		Thus, from the exact sequences \eqref{eq3} and \eqref{eq4} we conclude that
		\begin{align*}
			\Br(U/\Gamma)\ \simeq\ \dfrac{H^2(\Gamma,\, \bb{C}^*)}{H}\ \simeq\ \dfrac{H^2(\Gamma,\,
				\bb{C}^*)}{H'}\ \simeq\ \ker\left(\Br(V/\Gamma)\ \longrightarrow\ \Br(V)\right).
		\end{align*}
		This proves the proposition.
	\end{proof}
	
	The following is the main result of this section.
	
	\begin{theorem}\label{thm:brauer-group-concentrated-weights}
		Fix an even positive integer $r$, a line bundle $L$ on $X$, a set of parabolic points $S$ on $X$, a system of
		multiplicities $\textbf{m}$ of symmetric type, and a concentrated system of weights
		$\boldsymbol{\alpha}$ (Definition \ref{def:concentrated-weights}). Let
		$\left(\mc{N}^{\boldsymbol{m,\alpha},d}_{L(S)}\right)^{sm}$ denote the smooth locus of
		$\mc{N}^{\boldsymbol{m,\alpha},d}_{L(S)}$ (see Definition \ref{def:psp-moduli}).
		Then the Brauer group of $\left(\mc{N}^{\boldsymbol{m,\alpha},d}_{L(S)}\right)^{sm}$
		has the following description:
		\begin{enumerate}[(1)]
			\item \label{f1} If $d\,=\,0$ (equivalently, $\deg(L)$ is even),\, $\frac{r}{2}\,\geq\, 3$ is odd
			and $m^i_{_{p}}\,=\, 1$ for some $p\,\in\, S$ and $i$,
			\begin{align}
				\Br\left(\left(\mc{N}^{\boldsymbol{m,\alpha},d}_{L(S)}\right)^{sm}\right)
				\,\ \stackrel{\simeq}{\longrightarrow}\ \, \dfrac{H^2(\Gamma,\,\bb{C}^*)}{\frac{\bb{Z}}{2\bb{Z}}}.
			\end{align}
			
			\item \label{f2} If $d\,=\,0$ (equivalently, $\deg(L)$ is even),\, $\frac{r}{2}\,\geq\, 3$
			is even and $m^i_{_{p}}\,=\,1$ for some $p\,\in\, S$ and $i$,
			\begin{align}
				\Br\left(\left(\mc{N}^{\boldsymbol{m,\alpha},d}_{L(S)}\right)^{sm}\right)\ \,\stackrel{\simeq}{\longrightarrow}
				\ \, H^2(\Gamma,\,\bb{C}^*).
			\end{align}
			
			\item \label{f3} If $d\,=\,1$ (equivalently $\deg(L)$ is odd),\, $\frac{r}{2}\,\geq\, 3$ is even and
			$m^i_{_{p}}\,=\, 1$ for some $p\,\in\, S$ and $i$,
			\begin{align}
				\Br\left(\left(\mc{N}^{\boldsymbol{m,\alpha},d}_{L(S)}\right)^{sm}\right)\ \,\stackrel{\simeq}{\longrightarrow}
				\ \, H^2(\Gamma,\,\bb{C}^*).
			\end{align}
			
			\item \label{f4} If $d\,=\,1$ (equivalently, $\deg(L)$ is odd) and $\frac{r}{2}\,\geq\, 3$ is odd,
			\begin{align}
				\Br\left(\left(\mc{N}^{\boldsymbol{m,\alpha},d}_{L(S)}\right)^{sm}\right)\ \,\stackrel{\simeq}{\longrightarrow}
				\ \, H^2(\Gamma,\,\bb{C}^*).
			\end{align}
		\end{enumerate}
	\end{theorem}
	
	\begin{proof}
		By Proposition \ref{prop:brauer-group-U}, under any of the conditions
		\eqref{c1}, \eqref{c2} or \eqref{c3} we have
		\begin{align}\label{eqn:iso}
			\Br\left(U/\Gamma\right)\ \stackrel{\simeq}{\longrightarrow}\ \ker\left(\Br(V/\Gamma)\,\longrightarrow\,
			\Br(V)\right).
		\end{align}
		Now, $\ker\left(\Br(V/\Gamma)\longrightarrow\Br(V)\right)$ can be computed using
		\cite[Proposition 8.2]{BHol13}. To be more precise, the following exact sequence can be obtained using
		\cite[(8.1)]{BHol13}:
		\begin{align}\label{eqn:exact-seq-1}
			0\ \longrightarrow \ \bb{Z}/m\bb{Z}\ \longrightarrow\ H^2(\Gamma,\ \bb{C}^*)\ \longrightarrow\ 
			\Br\left(\left(\mc{N}^{d}_L\right)^{sm}\right)\ \longrightarrow\ \Br\left(\mc{M}^{rs}_L\right)\ \longrightarrow\ 0,
		\end{align}
		where $m$ is the smallest power of certain generating line bundle on the affine Grassmannian which
		descends (see \cite[Proposition 8.2]{BHol13} for details). Let us also recall the isomorphisms of Brauer groups obtained from Lemma \ref{lem:picard-Brauer-V} and equation \eqref{eqn:codimension-estimate-V-mod-Gamma}, namely
		$$\Br(V)\ \stackrel{\simeq}{\longrightarrow} \ \Br(\mc{M}^{rs}_L)\ \ \text{and}\ \ \Br(V/\Gamma)\ \stackrel{\simeq}{\longrightarrow}\ \Br\left(\left(\mc{N}^{d}_L\right)^{sm}\right) .$$ 
		Combining these isomorphisms with the exact sequence \eqref{eqn:exact-seq-1} enables us to conclude that
		$$\left(\ker\left(\Br(V/\Gamma)\,\longrightarrow\,
		\Br(V)\right)\right)\ \stackrel{\simeq}{\longrightarrow}\ \dfrac{H^2(\Gamma,\ \bb{C}^*)}{\frac{\bb{Z}}{m\bb{Z}}}\ \ (m \ \text{is as in}\ \eqref{eqn:exact-seq-1}). $$ 
		
		Now, if $d\,=\,0$ and $\frac{r}{2}\geq 3$ is odd, it follows from \cite{BHol13} that $m\,=\,2$ in
		\eqref{eqn:exact-seq-1}. If one moreover assumes that $m_{_{p}}^i\,=\,1$ for some $p\,\in\, S$ and $i$, using
		Proposition \ref{prop:brauer-group-U} the following is obtained:
		\begin{align*}
			\Br\left(\left(\mc{N}^{\boldsymbol{m,\alpha},d}_{L(S)}\right)^{sm}\right)\ \underset{\eqref{eqn:brauer-group-U}}{\simeq}\ \Br(U/\Gamma)\ \underset{\eqref{eqn:iso}}{\simeq}\ 
			\ker\left(\Br(V/\Gamma)\ \longrightarrow\ \Br(V)\right)\ \simeq\ \dfrac{H^2(\Gamma,\ \bb{C}^*)}{\frac{\bb{Z}}{2\bb{Z}}}
			\ .
		\end{align*}
		This proves the case \eqref{f1} of the theorem.
		
		Regarding the remaining cases \eqref{f2}, \eqref{f3} and \eqref{f4}, it follows from
		\cite{BHol13} that $m\,=\,1$ in \eqref{eqn:exact-seq-1}. Thus, in each of the remaining cases \eqref{f2},
		\eqref{f3} and \eqref{f4}, using Proposition \ref{prop:brauer-group-U} one obtains the following:
		\begin{align*}
			\Br\left(\left(\mc{N}^{\boldsymbol{m,\alpha},d}_{L(S)}\right)^{sm}\right)\ \underset{\eqref{eqn:brauer-group-U}}{\simeq}\ \Br(U/\Gamma)\ \underset{\eqref{eqn:iso}}{\simeq}\ 
			\ker\left(\Br(V/\Gamma)\ \longrightarrow\ \Br(V)\right)\ \simeq\ H^2(\Gamma,\ \bb{C}^*)
			\ .
		\end{align*}
		This proves the theorem.
	\end{proof}
	
	\section{Brauer groups for arbitrary systems of weights}
	
	The previous section dealt with concentrated systems of weights (Definition \ref{def:concentrated-weights}). We will now address the 
	situation where the system of weights $\boldsymbol{\alpha}$ need not be concentrated. In order to do so, we 
	first make a few remarks regarding the construction of the parabolic symplectic moduli space 
	$\mc{M}^{\boldsymbol{m,\alpha}}_{L(S)}$.
	
	Let $G$ be a connected reductive affine algebraic group over $\mathbb C$ acting on a projective variety $Y$. 
	In order to construct a geometric invariant theoretic quotient of $Y$ under the action of $G$, one fixes an
	ample 
	$G$--equivariant line bundle on $Y$. A natural question: How the quotient changes as the
	$G$--equivariant line bundle changes?
	Various authors including Boden--Hu, Dolgachev--Hu, Thaddeus and others 
	have studied this question. There are notions of 
	chambers and walls in the the $G$--ample cone in the N\'eron--Severi group of $G$--linearized line 
	bundles on $Y$ (\cite[Definition 0.2.1]{DolHu98}, \cite{Th96}); the geometric invariant theoretic quotient
	does not change as long as the line bundle remains in the interior of a chamber.
	
	The moduli space $\mc{M}^{\boldsymbol{m,\alpha}}_{L(S)}$ has been constructed and studied in 
	\cite{WaWe24} under the assumption on the system of weights and multiplicities that they are of 
	symmetric type (see \cite[Definition 2.2]{WaWe24}).
	
	Fix a set of parabolic points $S$ and also a system of multiplicities $\boldsymbol{m}$ at these points. Consider a 
	system of weights which is compatible with $\boldsymbol{m}$. If the system of weights consists of 
	\textit{rational} numbers, then such a choice of weights amounts to choosing a polarization on a certain 
	product of flag varieties for taking the GIT quotient by a suitable special linear group (see \cite[\S 
	~3]{WaWe24}). Thus, the set of all possible system of weights of symmetric type which are compatible with 
	$\boldsymbol{m}$ correspond to elements in the cone of ample linearized line bundles mentioned above (see 
	\cite{DolHu98,Th96}). By the virtue of variation of GIT principles, this cone is separated by finitely 
	many hyperplanes called \textit{walls}, and the connected components of these hyperplane complements are 
	known as \textit{chambers}. The moduli space remains unchanged as long as the system of weights vary 
	inside a chamber. We shall call a system of weights as \textit{generic} if it is contained in a chamber. 
	Now, since the collection of \textit{concentrated} systems of weights is clearly an open subset in this 
	cone, and the intersections of walls are of codimension one, clearly there exists a \textit{concentrated} 
	system of weights inside the cone which is not contained in any wall. We thus conclude that there exists a 
	\textit{generic} concentrated system of weights.
	
	\begin{proposition}\label{prop:variation-of-weights}
		Fix a system of multiplicities $\boldsymbol{m}$, and let $\boldsymbol{\alpha}$ and $\boldsymbol{\beta}$ be
		two systems of weights compatible with $\boldsymbol{m}$ in adjacent chambers in the ample cone which are
		separated by a single wall. Let $\mc{M}^{\boldsymbol{m,\alpha}}_{L(S)}$ and
		$\mc{M}^{\boldsymbol{m,\beta}}_{L(S)}$ denote the corresponding moduli spaces of semistable parabolic
		symplectic vector bundles. Then 
		$$\Br\left(\left(\mc{N}^{\boldsymbol{m,\alpha},d}_{L(S)}\right)^{sm}\right)\ \ \simeq\ \
		\Br\left(\left(\mc{N}^{\boldsymbol{m,\beta},d}_{L(S)}\right)^{sm}\right).$$
	\end{proposition}
	
	\begin{proof}
		Let $U_{\boldsymbol{\alpha}}$ denote the open subset of $\mc{M}^{\boldsymbol{m,\alpha}}_{L(S)}$ consisting
		of those stable parabolic symplectic vector bundles of quasi--parabolic type $\boldsymbol{m}$ that are
		both $\boldsymbol{\alpha}$--stable as well as $\boldsymbol{\beta}$--stable. Similarly, let
		$U_{\boldsymbol{\beta}}$ denote the open subset of $\mc{M}^{\boldsymbol{m,\beta}}_{L(S)}$ consisting of
		those stable parabolic symplectic vector bundles of quasi--parabolic type $\boldsymbol{m}$ that are both
		$\boldsymbol{\alpha}$--stable as well as $\boldsymbol{\beta}$--stable. By \cite[Theorem 3.5]{Th96}, there
		exists a birational morphism between these two moduli, and moreover, this birational morphism restricts
		to an isomorphism between $U_{\boldsymbol{\alpha}}$ and $U_{\boldsymbol{\beta}}$, which we denote by $\textbf g$:
		\begin{equation}\label{eg}
			{\textbf g}\ :\ U_{\boldsymbol{\alpha}}\ \, \overset{\simeq}{\longrightarrow}\
			\, U_{\boldsymbol{\beta}}.
		\end{equation}
		This isomorphism is given simply by interchanging the weights between
		$\boldsymbol{\alpha}$ and $\boldsymbol{\beta}$, keeping the underlying quasi-parabolic vector
		bundle unchanged. Moreover, the complement of $U_{\boldsymbol{\alpha}}$ in $\mc{M}^{\boldsymbol{m,\alpha}}_{L(S)}$
		is of codimension at least 2, and the same holds for the complement of
		$U_{\boldsymbol{\beta}}$ in $\mc{M}^{\boldsymbol{m,\beta}}_{L(S)}$.
		
		Now, assume that $\boldsymbol{\alpha}$ is concentrated in the sense of Definition \ref{def:concentrated-weights}.
		Define $$U^{sm}_{\boldsymbol{\alpha}}\ \, :=\ \, U_{\boldsymbol{\alpha}}\cap
		\left(\mc{M}^{\boldsymbol{m,\alpha}}_{L(S)}\right)^{sm}.$$
		It follows from a straightforward dimension comparison that the complement of $U^{sm}_{\boldsymbol{\alpha}}$
		in $\left(\mc{M}^{\boldsymbol{m,\alpha}}_{L(S)}\right)^{sm}$ is of codimension at least 2 as well. 
		
		As before, $\Gamma$ denotes the group of $2$--torsion line bundles on $X$. It is clear from their description
		that both $U_{\boldsymbol{\alpha}}$ and $U_{\boldsymbol{\beta}}$ are $\Gamma$--invariant open subsets
		of the moduli spaces, and the isomorphism $\textbf g$ in \eqref{eg} is $\Gamma$--equivariant as well.
		Now, as we have seen in the proof of Theorem \ref{thm:brauer-group-concentrated-weights}, there is a
		Zariski open subset $U\,\subset\, \left(\mc{M}^{\boldsymbol{m,\alpha}}_{L(S)}\right)^{sm}$ whose complement
		has codimension at least 2, and moreover $\Gamma$ acts freely on $U$. Let $$U'\,\ :=\,\ 
		U\cap U^{sm}_{\boldsymbol{\alpha}}.$$
		We can write
		$$\left(\mc{M}^{\boldsymbol{m,\alpha}}_{L(S)}\right)^{sm}\setminus U'\ =\
		\left(\left(\mc{M}^{\boldsymbol{m,\alpha}}_{L(S)}\right)^{sm}\setminus U\right) \cup
		\left(\left(\mc{M}^{\boldsymbol{m,\alpha}}_{L(S)}\right)^{sm}\setminus U_{\boldsymbol{\alpha}}\right).$$
		Since both $\left(\left(\mc{M}^{\boldsymbol{m,\alpha}}_{L(S)}\right)^{sm}\setminus U\right)$, as well as $\left(\left(\mc{M}^{\boldsymbol{m,\alpha}}_{L(S)}\right)^{sm}\setminus U_{\boldsymbol{\alpha}}\right)$, are closed subsets of codimension at least 2
		in $\left(\mc{M}^{\boldsymbol{m,\alpha}}_{L(S)}\right)^{sm}$, it follows that the complement of $U'$ in
		$\left(\mc{M}^{\boldsymbol{m,\alpha}}_{L(S)}\right)^{sm}$ is of codimension at least 2 as well. Moreover,
		being an intersection of two $\Gamma$--invariant open subsets, $U'$ is also $\Gamma$--invariant.
		Since $\Gamma$ acts freely on $U$, it follows that the action of $\Gamma$ on $U'$ is free. 
		
		For the isomorphism $\textbf{g}$ in \eqref{eg},
		the image $\textbf{g}(U')$ is again a $\Gamma$--invariant Zariski open subset of
		$\left(\mc{M}^{\boldsymbol{m,\beta}}_{L(S)}\right)^{sm}$ whose complement is of codimension at least $2$
		on which $\Gamma$ acts freely. Upon taking quotient by the $\Gamma$--action, the $\Gamma$--equivariant
		isomorphism
		$$\textbf{g}\big\vert_{U'} \ :\ U'\ \stackrel{\simeq}{\longrightarrow} \ \textbf{g}(U')$$ descends to an isomorphism between $U'/\Gamma$
		and $\textbf{g}(U')/\Gamma$. Clearly,
		$U'/\Gamma \,\subset\, \left(\mc{N}^{\boldsymbol{m,\alpha},d}_{L(S)}\right)^{sm}$ and $\textbf{g}(U')/\Gamma\,
		\subset\,\left(\mc{N}^{\boldsymbol{m,\beta},d}_{L(S)}\right)^{sm}$ are Zariski open subsets whose complements
		are of codimension at least $2$. Consequently, we have
		\begin{align}
			\Br\left(\left(\mc{N}^{\boldsymbol{m,\alpha},d}_{L(S)}\right)^{sm}\right)\, \ \simeq\,\ \Br\left(U'/\Gamma\right)\
			\,\simeq\, \ \Br\left(\textbf{g}(U')/\Gamma\right)\,\ \simeq\,\
			\Br\left(\left(\mc{N}^{\boldsymbol{m,\beta},d}_{L(S)}\right)^{sm}\right).
		\end{align}
		Next, since there are only finitely many chambers and walls, one can arrange the collection
		of chambers in a sequence, say $C_1,\ \cdots,\ C_N$, where $C_1$ contains a concentrated system of
		weights (Definition \ref{def:concentrated-weights}), and for each $1\,\leq\, i\,\leq\, N$, the chambers $C_i$
		and $C_{i+1}$ are separated by a single wall. Then, one can inductively go from $C_1$ to $C_2$, then from $C_2$
		to $C_3$ and so on. This proves the proposition.
	\end{proof}
	
	\begin{corollary}\label{cor:arbitrary-weight}
		Under any of the conditions $(1)-(4)$ in Theorem \ref{thm:brauer-group-concentrated-weights}, the conclusion of Theorem \ref{thm:brauer-group-concentrated-weights} remains valid for any generic system of weights in the ample cone.
	\end{corollary}
	
	\begin{proof}
		This follows immediately from Proposition \ref{prop:variation-of-weights}.
	\end{proof}

\end{document}